\def\yes{\if00}

\def\iftwelvept{\yes}
\def\ifusepdf{\yes}
\def\ifpsfont{\yes}

\iftwelvept
\documentclass[leqno,12pt]{amsart}
\else
\documentclass[leqno]{amsart}
\fi
\usepackage{amssymb}
\usepackage{amscd}
\usepackage{latexsym}
\usepackage{verbatim}
\ifusepdf
\usepackage{hyperref}
\else\fi
\ifpsfont
\usepackage[T1]{fontenc}
\usepackage{times}
\else\fi

\iftwelvept
\else\fi

\theoremstyle{plain}
\newtheorem{Theorem}{Theorem}[section]
\newtheorem{Proposition}[Theorem]{Proposition}
\newtheorem{Lemma}[Theorem]{Lemma}
\newtheorem{Corollary}[Theorem]{Corollary}
\newtheorem{Claim}{Claim}[Theorem]

\theoremstyle{definition}
\newtheorem{Definition}[Theorem]{Definition}
\newtheorem{Remark}[Theorem]{Remark}

\renewcommand{\theTheorem}{\arabic{section}.\arabic{Theorem}}
\renewcommand{\theClaim}{\arabic{section}.\arabic{Theorem}.\arabic{Claim}}
\renewcommand{\theequation}{\arabic{section}.\arabic{Theorem}.\arabic{Claim}}

\def\rom{\textup}
\newcommand{\ZZ}{{\mathbb{Z}}}
\newcommand{\QQ}{{\mathbb{Q}}}
\newcommand{\RR}{{\mathbb{R}}}
\newcommand{\CC}{{\mathbb{C}}}
\newcommand{\DD}{{\mathbb{D}}}
\newcommand{\PP}{{\mathbb{P}}}
\newcommand{\NN}{{\mathbb{N}}}
\newcommand{\OO}{{\mathcal{O}}}
\newcommand{\XX}{{\mathcal{X}}}
\newcommand{\LL}{{\mathcal{L}}}

\newcommand{\ord}{\operatorname{ord}}

\newcommand{\Coker}{\operatorname{Coker}}
\newcommand{\Image}{\operatorname{Image}}
\newcommand{\length}{\operatorname{length}}

\newcommand{\Supp}{\operatorname{Supp}}

\newcommand{\Pic}{\operatorname{Pic}}
\newcommand{\aPic}{\widehat{\operatorname{Pic}}}
\newcommand{\Spec}{\operatorname{Spec}}

\newcommand{\aChow}{\widehat{\operatorname{CH}}}
\newcommand{\aCycle}{\widehat{Z}}
\newcommand{\Proj}{\operatorname{Proj}}
\newcommand{\zero}{\operatorname{div}}
\newcommand{\Proof}{{\sl Proof.}\quad}
\newcommand{\adeg}{\widehat{\operatorname{deg}}}
\newcommand{\trdeg}{\operatorname{tr.deg}}
\newcommand{\rank}{\operatorname{rk}}
\newcommand{\acherncl}{\widehat{{c}}}
\newcommand{\normabb}{\Vert\!\cdot\!\Vert}
\newcommand{\Bs}{\operatorname{Bs}}
\newcommand{\SBs}{\operatorname{SBs}}
\newcommand{\QED}{{\unskip\nobreak\hfil\penalty50\quad\null\nobreak\hfil
{$\Box$}\parfillskip0pt\finalhyphendemerits0\par\medskip}}
\newcommand{\rest}[2]{\left.{#1}\right\vert_{{#2}}}
\begin{document}

%%%%%%%%%%%
%% Title %%
%%%%%%%%%%%
\title[Arithmetic height functions over finitely generated fields]%
{Arithmetic height functions \\ over finitely generated fields}
\author{Atsushi Moriwaki}
\address{Department of Mathematics, Faculty of Science,
Kyoto University, Kyoto, 606-8502, Japan}
\email{moriwaki@kusm.kyoto-u.ac.jp}
%\date{\DateTime, (\Version)}
\date{1/June/1999, 1:50PM (JP), (Version 2.0)}
\begin{abstract}
In this paper, we propose a new height function
for a variety defined over a finitely generated field over $\QQ$.
For this height function, we will prove Northcott's theorem and
Bogomolov's conjecture, so that we can recover
the original Raynaud's theorem (Manin-Mumford's conjecture).
\end{abstract}

%%%%%
% 1st Draft
%%%%%
%\vskip -1cm
%\hfill\fbox{{\large\bf 1st draft}} %
%\footnote{I welcome any comments and suggestions.}\par
%\vskip .5cm
%%%%%

\maketitle

\tableofcontents

\renewcommand{\theTheorem}{\Alph{Theorem}}

\section*{Introduction}
Let $K$ be a finitely generated field over $\QQ$,
and $d$ the transcendence degree of $K$ over $\QQ$.
If $d=1$,
then there is a smooth projective curve $C$ over
a number field such that the function field of $C$ is $K$.
Using non-archimedean valuations arising from
points of $C$, we can define a geometric height function
\[
h^{geom} : \PP^n(\overline{K}) \to \RR.
\]
It is well know that this height function can be given in terms of
the usual intersection theory, so that
it is rather easy to handle it.
However,
in contract with
height functions over
number fields, it does not
reflect the exact state of points on $\PP^n(\overline{K})$.
For example, Northcott's theorem does not hold for
the geometric height function in general.
A reason for this, we can consider, is that
$h^{geom}$ does not take care of data coming from
the constant field.
In this paper, we would like to propose
a new kind of height functions for finitely generated fields over $\QQ$,
and unify them with the geometric height functions.

A key idea to get a new height function 
is to fix a polarization $\overline{B} = (B; \overline{H}_1, \ldots, \overline{H}_d)$ of $K$,
namely, a collection of a normal projective arithmetic variety $B$
whose function field is $K$, and nef $C^{\infty}$-hermitian line bundles
$\overline{H}_1, \ldots, \overline{H}_d$
on $B$. Here a $C^{\infty}$-hermitian line bundle $\overline{H}$
is said to be nef if $c_1(\overline{H})$ is semipositive, and
$\adeg\left( \rest{\overline{H}}{\Gamma} \right) \geq 0$
for any one-dimensional integral subschemes $\Gamma$ of $B$.
Once we fix the polarization $\overline{B}$ of $K$,
then we can define a height function 
\[
h^{\overline{B}}_K : \PP^n(K) \to \RR
\]
associated with $\overline{B}$ to be
\begin{multline*}
h^{\overline{B}}_K(\phi_0, \ldots, \phi_n) =
\sum_{\Gamma} \max_{i} \{ - \ord_{\Gamma}(\phi_i) \}
\adeg \left( \acherncl_1\left(\rest{\overline{H}_1}{\Gamma}\right) \cdots
\acherncl_1\left(\rest{\overline{H}_d}{\Gamma}\right) \right) \\
+ \int_{B(\CC)} \log\left( \max_i \{ \vert \phi_i \vert \} \right)
c_1(\overline{H}_1) \wedge \cdots \wedge c_1(\overline{H}_d),
\end{multline*}
where $\Gamma$ runs over all prime divisors on $B$.
Moreover, we can easily see that $h^{\overline{B}}_K$ extends to
\[
h^{\overline{B}} : \PP^n(\overline{K}) \to \RR.
\]
For example, if $d=1$ and $\overline{H}_1$ is given 
by the infinite fibers of $B$,
then $h^{\overline{B}}$ is nothing more than $h^{geom}$ up to
the multiplication of a positive constant.
Moreover, note that if $d=0$, then $h^{\overline{B}}$ is the
usual height function over a number field.

Further, we can give these height functions in terms of Arakelov
intersection theory.
Let $X$ be a projective variety over $K$, and $L$ a line bundle
on $X$. Let us take a model $(\XX, \overline{\LL})$ of $(X, L)$,
namely, $\XX$ is a projective arithmetic variety over $B$ and
$\overline{\LL}$ is a hermitian line bundle on $\XX$ with
$\XX_K = X$ and $\LL_K = L$.
For a point $P \in X(\overline{K})$,
we denote by $\Delta_P$ the closure of the image of
$\Spec(\overline{K}) \overset{P}{\longrightarrow} X \hookrightarrow \XX$.
Then, we define
\[
h^{\overline{B}}_{(\XX, \overline{\LL})} : X(\overline{K}) \to \RR
\]
to be
\[
h^{\overline{B}}_{(\XX, \overline{\LL})}(P) =
\frac{1}{[K(P):K]}\adeg \left(
\acherncl_1(\rest{\overline{\LL}}{\Delta_P}) \cdot
\acherncl_1(\rest{f^* (\overline{H}_1)}{\Delta_P}) \cdots
\acherncl_1(\rest{f^* (\overline{H}_d)}{\Delta_P})
\right),
\]
where $f : \XX \to B$ is the canonical morphism.
We can see that $h^{\overline{B}}_{(\XX, \LL)}$ modulo
the set of bounded functions on $X(\overline{K})$
does not depend on the choice
of the model $(\XX, \overline{\LL})$ of $(X, L)$
(cf. Corollary~\ref{cor:height:well:def:const}),
so that we may denote $h^{\overline{B}}_{(\XX, \overline{\LL})}$
by $h^{\overline{B}}_L$.

Since our height functions include the geometric height functions,
Northcott's theorem does not hold in general.
However, if the polarization $\overline{B}$ is big,
we can expect a certain kind of affirmative answers.
For this reason,
we introduce the following notation.
If $(H_i)_{\QQ}$'s are big on $B_{\QQ}$ and
there are positive numbers $n_1, \ldots, n_d$ such that
$H^0(B, H_i^{\otimes n_i})$ has a strictly small section for each $i$,
then $h^{\overline{B}}_L$ is called an
{\em arithmetic height function} and is denoted by
$h^{arith}_L$ for simplicity.
Then, we have the following Northcott's theorem
for the arithmetic height function.

\begin{Theorem}[cf. Theorem~\ref{thm:northcott:thm:fun:field}]
\label{thm:northcott:thm:fun:field:intro}
If $L$ is ample, then,
for any numbers $M$ and any positive integers $e$, the set
\[
\{ P \in X(\overline{K}) \mid h^{arith}_L(P) \leq M,
\quad [K(P) : K] \leq e \} 
\]
is finite.
\end{Theorem}

Now let $A$ be an abelian variety over $K$, and
$L$ a symmetric ample line bundle on $A$.
Then, in the same way as the usual height theory,
we can assign the canonical height (N\'{e}ron-Tate height)
$\hat{h}^{arith}_L$ to $(A, L)$.
Then, $\hat{h}^{arith}_L(x) \geq 0$ for all $x \in A(\overline{K})$,
and $\hat{h}^{arith}_L(x) = 0$ if and only if $x$ is a torsion
point as a corollary of Theorem~\ref{thm:northcott:thm:fun:field:intro}
(cf. Proposition~\ref{prop:positivity:canonical:height}).
Moreover, in terms of $\hat{h}^{arith}_L$,
we have the following solution of 
Bogomolov's conjecture over $K$, which is a generalization of
results due to Ullmo \cite{UlPos} and Zhang \cite{ZhEqui}.

\begin{Theorem}[cf. Theorem~\ref{thm:bogomolov:conj:fun}]
\label{thm:bogomolov:conj:fun:intro}
Let $X$ be a subvariety of $A_{\overline{K}}$.
If the set
\[
 \{ P \in X(\overline{K}) \mid \hat{h}^{arith}_L(P) \leq \epsilon \}
\]
is Zariski dense in $X$ for any positive numbers $\epsilon$,
then $X$ is a translation of an abelian subvariety of $A_{\overline{K}}$ 
by a torsion point.
\end{Theorem}

As corollary, we can recover
the original Raynaud's theorem (\cite{Ray1} and \cite{Ray2})
conjectured by Manin and Mumford.

%%
%\begin{Corollary}[cf. Corollary~\ref{cor:raynaud:thm}]
%\label{cor:raynaud:thm:intro}
%Let $A$ be an abelian variety over the complex number field $\CC$, and
%$X$ a subvariety of $A$. If $X(\CC) \cap A(\CC)_{tor}$ is Zariski dense
%in $X$, then $X$ is a translation of an abelian subvariety of $A$
%by a torsion point.
%\end{Corollary}
%%
\begin{Corollary}[cf. Corollary~\ref{cor:raynaud:thm}]
\label{cor:raynaud:thm:intro}
Let $A$ be an abelian variety over the complex number field $\CC$, and
$Z$ a reduced subscheme of $A$.
Then, every irreducible component of
the Zariski closure of $Z(\CC) \cap A(\CC)_{tor}$ in $A$ is 
a translation of an abelian subvariety of $A$
by a torsion point.
Consequently, there are finitely many abelian subvarieties
$B_1, \ldots, B_n$ of $A$ and torsion points $b_1, \ldots, b_n$ of $A(\CC)$
such that
\[
\overline{Z(\CC) \cap A(\CC)_{tor}} = \bigcup_{i=1}^n (B_i(\CC) + b_i)
\quad\text{and}\quad
Z(\CC) \cap A(\CC)_{tor} = \bigcup_{i=1}^n (B_i(\CC)_{tor} + b_i).
\]
\end{Corollary}

Finally, we would like to express gratitude to
Dr. Kawaguchi, Prof. Poonen, Prof. Szpiro, Prof. Ullmo, and Prof. Zhang for
their helpful conversations.
The author also thanks Prof. Silverman for his nice comments.

\renewcommand{\theTheorem}{\arabic{section}.\arabic{Theorem}}
\renewcommand{\theClaim}{\arabic{section}.\arabic{Theorem}.\arabic{Claim}}
\renewcommand{\theequation}{\arabic{section}.\arabic{Theorem}.\arabic{Claim}}

\section{Arakelov intersection theory}
In this paper, an arithmetic variety means a
flat and
quasi-projective integral scheme over $\ZZ$.
Moreover, we say an arithmetic variety is {\em generically smooth}
if it is smooth over $\QQ$.
For basic materials of Arakelov intersection theory,
we refer to \cite{GSArInt} and \cite{SoAr}.

Let $X$ be a generically smooth arithmetic variety.
According to \cite{KMSemi}, 
a pair $(Z, g)$ is called an {\em arithmetic cycle of codimension $p$}
(resp. {\em arithmetic $D$-cycle of codimension $p$})
if $Z$ is a cycle of codimension $p$ on $X$, and
$g$ is a Green current for $Z(\CC)$ (resp.
$g$ is a current of type $(p-1,p-1)$ on $X(\CC)$).
The set of all arithmetic cycles (resp. $D$-cycles) of codimension $p$
is denoted by $\aCycle^p(X)$ (resp. $\aCycle_D^p(X)$).
Let $\widehat{R}^p(X)$ be the subgroup of $\aCycle^p(X)$ generated
by the following elements:
\begin{enumerate}
\renewcommand{\labelenumi}{(\roman{enumi})}
\item 
$((f), - [\log |f|^2])$, 
where $f$ is a rational function on some
subvariety $Y$ of codimension $p-1$ and $[\log |f|^2]$ 
is the current defined by
\[
[\log |f|^2](\gamma) = 
        \int_{Y(\CC)} (\log |f|^2)\gamma.
\]

\item
$(0, \partial(\alpha) + \bar{\partial}(\beta))$,
where $\alpha$ and $\beta$ are currents of type $(p-2,p-1)$
and $(p-1,p-2)$ respectively.
\end{enumerate}
Here we define
\[
\aChow^p(X) = \aCycle^{p}(X)/\widehat{R}^p(X),
\quad\text{and}\quad
\aChow_D^p(X) = \aCycle_D^{p}(X)/\widehat{R}^p(X).
\]

Let $\overline{L} = (L, \normabb)$ be a 
$C^{\infty}$-hermitian line bundle on $X$.
We define a homomorphism
\renewcommand{\theequation}{\arabic{section}.\arabic{Theorem}}
\addtocounter{Theorem}{1}
\begin{equation}
\label{eqn:scalar:product:hermitian:line:bundle}
\acherncl_1(\overline{L}) \cdot\ : \aChow_D^p(X) \to \aChow_D^{p+1}(X)
\end{equation}
\renewcommand{\theequation}{\arabic{section}.\arabic{Theorem}.\arabic{Claim}}%
in the following way.
Let $(Z, g)$ be an element of $\aCycle_D^p(X)$.
We assume that $Z$ is integral.
Then, taking a rational section $s$ of $\rest{L}{Z}$,
we consider an arithmetic $D$-cycle
\[
(\zero(s), \left[- \log\left( \Vert s \Vert_Z^2 \right) \right] + 
c_1(\overline{L}) \wedge g),
\]
where
$\left[ -\log \left( \Vert s \Vert_Z^2 \right) \right]$
is a current given by
$\phi \mapsto -\int_{Z(\CC)} \log \left( \Vert s \Vert_Z^2 \right)
\phi$.
The class of the above cycle in $\aChow^{p+1}_D(X)$ does not depend on
the choice of the rational section $s$.
Thus, by linearity, we have a homomorphism
\[
\acherncl_1(\overline{L}) \cdot\ : \aCycle_D^p(X) \to \aChow_D^{p+1}(X).
\]
On the other hand, it is well known that
$\acherncl_1(\overline{L}) \cdot \widehat{R}^p(X) 
\subset \widehat{R}^{p+1}(X)$.
Thus, we obtain our desired homomorphism 
\eqref{eqn:scalar:product:hermitian:line:bundle}.  

Now let $\overline{M} = (M, \normabb)$ be a continuous hermitian line
bundle on $X$, namely, $\normabb$ is a continuous metric.
Then, $\acherncl_1(\overline{M})$ is defined by the class of
$(\zero(s), -\log \Vert s \Vert^2)$ in $\aChow_D^1(X)$,
where $s$ is a non-zero rational section of $M$.
This is actually well defined because the class does not depend on the
choice of the rational section $s$.
Then, using the scalar product
\eqref{eqn:scalar:product:hermitian:line:bundle},
for $C^{\infty}$-hermitian line bundles $\overline{L}_1,
\ldots, \overline{L}_d$ on $X$, we can define
\[
\acherncl_1(\overline{L}_1) \cdots \acherncl_1(\overline{L}_d)
\cdot \acherncl_1(\overline{M}) 
\in \aChow^{d+1}_D(X),
\]
where $d = \dim X_{\QQ}$. 
In particular, if $X$ is projective, then
we get the intersection number
\[
 \adeg \left(
\acherncl_1(\overline{L}_1) \cdots \acherncl_1(\overline{L}_d)
\cdot \acherncl_1(\overline{M}) 
 \right),
\]
where $\adeg :  \aChow^{d+1}_D(X) \to \RR$ is given by
\[
\adeg\left( \sum_P n_P P, T \right) = 
\sum_P n_P \log\#(\kappa(P)) + \frac{1}{2} \int_{X(\CC)} T.
\]

\medskip
Next, let us consider the push-forward of cycles.
Let $f : X \to Y$ be a projective morphism of 
generically smooth arithmetic varieties.
Then, $f_* : \aCycle_D^p(X) \to \aCycle_D^{p + \dim Y - \dim X}(Y)$
is defined by $f_*(Z, g) = (f_*(Z), f_*(g))$.
This induces
\[
f_* : \aChow_D^p(X) \to \aChow_D^{p + \dim Y - \dim X}(Y).
\]
Then, we have the following projection formula
(cf. \cite[Proposition~2.4.1]{KMSemi}).

\begin{Proposition}
\label{prop:projection:formular:general}
Let $\overline{L}$ be a $C^{\infty}$-hermitian line bundle on $Y$,
and $z$ an element of $\aChow_D^p(X)$.
Then,
\[
 f_*(\acherncl_1(f^*(\overline{L})) \cdot z) = 
 \acherncl_1(\overline{L}) \cdot f_*(z).
\]
\end{Proposition}

\Proof
For reader's convenience, we give the proof of it.
Let $(Z, g)$ be a representative of $z$, and
$\normabb$ the metric of $\overline{L}$.
Clearly, we may assume that $Z$ is reduced and irreducible.
We set $T = f(Z)$ and $\pi = \rest{f}{Z} : Z \to T$. 
Let $s$ be a non-zero rational section of $\rest{L}{T}$.
Then, $\pi^*(s)$ gives rise to a non-zero rational section of $\rest{f^*(L)}{Z} = 
\pi^* \left( \rest{L}{T} \right)$.
Thus, $\acherncl_1 (f^*(\overline{L})) \cdot z$ can be represented by
\[
\left( \zero(\pi^*(s)),\ \ %
\left[ -\log \pi^* \left(\Vert s \Vert_T^2  \right) \right]
+ c_1(f^*(\overline{L})) \wedge g \right).
\]
If we set
\[
\deg(\pi) =
\begin{cases}
0 & \text{if $\dim T < \dim Z$} \\
\deg(Z \to T) & \text{if $\dim T = \dim Z$,}
\end{cases}
\]
then
\begin{align*}
\int_{Z(\CC)}
-\log \pi^* \left( \Vert s \Vert_T^2 \right) f^*(\psi)
& = \int_{Z(\CC)} \pi^* \left( 
-\log \left( \Vert s \Vert_T^2 \right) \psi \right) \\
& = \deg(\pi) \int_{T(\CC)} -\log \left( \Vert s \Vert_T^2 \right) \psi
\end{align*}
for a compactly supported $C^{\infty}$-form $\psi$ on $Y(\CC)$.
Thus, we have
\[
f_* \left[ -\log \pi^* \left(\Vert s \Vert_T^2  \right) \right]
= \deg(\pi) \left[ -\log \left( \Vert s \Vert_T^2 \right) \right].
\]
Therefore,
\begin{align*}
f_*(\acherncl_1 (f^*(\overline{L})) \cdot z)
& =
\left( \deg(\pi) \zero(s),\ \ %
\deg(\pi) \left[ -\log \left( \Vert s \Vert_T^2 \right) \right] +
c_1(L, h) \wedge f_*(g) \right) \\
& = \acherncl_1 (\overline{L}) \cdot (\deg(\pi) T, f_*(g)) = \acherncl_1 (\overline{L}) \cdot f_*(z).
\end{align*}
Hence, we get our proposition.
\QED

\medskip
Finally, let us consider intersections on a general projective 
arithmetic variety.
Let $T$ be a reduced complex space, and $\overline{L} = (L, \normabb)$
a continuous hermitian line bundle on $T$.
According to \cite{ZhPos}, we say $\overline{L}$ is $C^{\infty}$
if, for any analytic morphisms $h : M \to T$
from any complex manifolds $M$ to $T$,
$h^*(\overline{L})$ is 
a $C^{\infty}$-hermitian line bundle on $M$.
In the same way, we can define $C^{\infty}$-functions on $T$.
From now on, we assume that $\overline{L}$ is $C^{\infty}$.
Then, $c_1(\overline{L})$ is a $C^{\infty}$-form on $T_{reg}$,
where $T_{reg}$ is the smooth locus of $T$.
We say $c_1(\overline{L})$ is {\em semipositive} if,
for any analytic morphisms $h : M \to T$
from any complex manifolds $M$ to $T$,
$c_1(h^*(\overline{L}))$ is a semipositive form on $M$.
Moreover, we say $c_1(\overline{L})$ is {\em positive} if,
for any real valued $C^{\infty}$-functions $f$ on $T$
with compact support,
there is a positive real number $\lambda_0$ such that
$c_1(\overline{L}) + \lambda dd^c(f)$ is semipositive
for all $\lambda$ with $|\lambda| \leq \lambda_0$.
Note that $dd^c(f)$ is the Chern form of
$(\OO_T, \exp(-f)|\cdot|)$.
It is easy to see that the above positivity of Chern forms coincides
with the usual positivity of them if $T$ is non-singular.
\begin{comment}
Moreover, we say $c_1(\overline{L})$ is {\em positive} if,
for any point of $t \in T$,
there are an open neighborhood $U$ of $t$,
a complex manifold $W$, and $C^{\infty}$ positive
$(1,1)$-form $\omega$ on $W$ such that
$U$ is an analytic closed subset of $W$ and
$c_1(\overline{L})$ is the restriction of $\omega$ to $U$.
We can see that
if $c_1(\overline{L})$ is positive, then it is semipositive.
\end{comment}

Let $X$ be a projective arithmetic variety of $d = \dim X_{\QQ}$.
Let $\mu : X' \to X$ be a generic resolution of singularities of $X$,
namely, $\mu : X' \to X$ is a birational morphism of 
projective arithmetic varieties such that $X'$ is generically smooth.
Note that a generic resolution of singularities exists
for any arithmetic varieties by using Hironaka's resolution
of singularities \cite{Hiro}.
Let $\overline{L}_1, \cdots, \overline{L}_d$
be $C^{\infty}$-hermitian line bundles on $X$, and
$\overline{M}$ a continuous hermitian line bundle on $X$.
Then, the intersection number
$\adeg \left( \acherncl_1(\mu^* \overline{L}_1) \cdots 
\acherncl_1(\mu^* \overline{L}_d) \cdot \acherncl_1(\mu^* \overline{M})
\right)$ does not depend on the choice of the generic resolution
of singularities $\mu : X' \to X$ by virtue of
Proposition~\ref{prop:projection:formular:general}, so that
we define
\[
 \adeg \left( \acherncl_1(\overline{L}_1) \cdots 
 \acherncl_1(\overline{L}_d) \cdot \acherncl_1(\overline{M})
 \right)
\] 
by
$\adeg \left( \acherncl_1(\mu^* \overline{L}_1) \cdots 
\acherncl_1(\mu^* \overline{L}_d) \cdot \acherncl_1(\mu^* \overline{M})
\right)$.
Then, we have the following proposition.

\begin{Proposition}
\label{prop:projection:formula:intersection:num}
Let $f : X \to Y$ be a morphism of projective arithmetic varieties
with $d = \dim X_{\QQ}$ and $n = \dim Y_{\QQ}$
\begin{enumerate}
\renewcommand{\labelenumi}{(\arabic{enumi})}
\item
Let $\overline{L}_1, \ldots, \overline{L}_r$
be $C^{\infty}$-hermitian line bundles on $X$, and
$\overline{M}_1, \ldots, \overline{M}_s$ $C^{\infty}$-hermitian
line bundles on $Y$. We assume that $r + s = d+1$.
Then,
\begin{multline*}
\adeg \left(
\acherncl_1(\overline{L}_1) \cdots \acherncl_1(\overline{L}_r) \cdot
\acherncl_1(f^* \overline{M}_1) \cdots \acherncl_1(f^* \overline{M}_s)
\right) \\
= \begin{cases}
0 & \text{if $s > n + 1$} \\
\deg((L_1)_{\eta} \cdots (L_r)_{\eta})\adeg \left(
\acherncl_1(\overline{M}_1) \cdots \acherncl_1(\overline{M}_s)
\right) & \text{if $s = n+1$},
\end{cases}
\end{multline*}
where the $\eta$ means the restriction to the generic fiber of $f$.
\item
We assume that $f$ is generically finite.
Let $\overline{L}_1, \ldots, \overline{L}_n$ be $C^{\infty}$-hermitian
line bundles on $Y$, and $\overline{M}$ a continuous hermitian line
bundle on $Y$.
Then,
\[
 \adeg \left(
 \acherncl_1(f^* \overline{L}_1) \cdots \acherncl_1(f^* \overline{L}_n)
 \cdot \acherncl_1(f^* \overline{M})
\right) = \deg(f)
\adeg \left(
 \acherncl_1(\overline{L}_1) \cdots \acherncl_1(\overline{L}_n)
 \cdot \acherncl_1(\overline{M})
\right).
\]
\end{enumerate}
\end{Proposition}

\Proof
Let us consider the following commutative diagram:
\[
 \begin{CD}
 X @<{\mu}<< X' \\
 @V{f}VV @VV{f'}V \\
 Y @<{\nu}<< Y' 
 \end{CD}
\]
where $\mu : X' \to X$ and $\nu : Y' \to Y$ are
generic resolutions of singularities of $X$ and $Y$ respectively.
Then, our assertions are consequences of 
Proposition~\ref{prop:projection:formular:general}.
\QED

%%%
\renewcommand{\theTheorem}{\arabic{section}.\arabic{Theorem}}
\renewcommand{\theClaim}{\arabic{section}.\arabic{Theorem}.\arabic{Claim}}
\renewcommand{\theequation}{\arabic{section}.\arabic{Theorem}.\arabic{Claim}}
%%%
\section{Arithmetically positive hermitian line bundles}
Let $B$ be a projective arithmetic variety with $d = \dim B_{\QQ}$, and
$\overline{H}$ a $C^{\infty}$-hermitian line bundle on $B$.
From now on, we introduce several kinds of positivity of $\overline{H}$.
First of all, we say $\overline{H}$ is {\em ample} if
$H$ is ample, $c_1(\overline{H})$ is semipositive on $B(\CC)$, and,
for a sufficiently large $n$,
$H^0(B, H^{\otimes n})$ is generated by
$\{ s \in H^0(B, H^{\otimes n}) \mid \Vert s \Vert_{\sup} < 1 \}$.
$\overline{H}$ is said to be {\em vertically nef}
if $H$ is relatively nef with respect to $B \to \Spec(\ZZ)$, and
$c_1(\overline{H})$ is semipositive on $B(\CC)$.
We say $\overline{H}$ is {\em horizontally nef} if,
for all one-dimensional integral closed subschemes $\Gamma$ 
flat over $\ZZ$,
$\adeg \left( \rest{\overline{H}}{\Gamma} \right) \geq 0$.
Moreover, if $\overline{H}$ is vertically nef and horizontally nef,
then we say $\overline{H}$ is {\em nef}.
Further, $\overline{H}$ is said to be {\em big} if 
$\rank H^0(B, H^{\otimes m}) = O(m^d)$, and there is a non-zero
section $s$ of $H^0(B, H^{\otimes n})$ with $\Vert s \Vert_{\sup} < 1$ for
some positive integer $n$.
It is easy to see that
if $\overline{H}$ is ample, then $\overline{H}$ is nef and big.
The following theorem due to Faltings, Gillet-Soul\'{e} and
Zhang is a very useful criterion for the bigness of $C^{\infty}$-hermitian
line bundles (cf. \cite[Theorem~1.4]{ZhPos}).

\begin{Theorem}
\label{thm:existense:small:sec}
Let $B$ be a projective arithmetic variety with $d = \dim B_{\QQ}$, and
$\overline{L}$ a $C^{\infty}$-hermitian line bundle on $B$.
If $\overline{L}$ is vertically nef, 
$L_{\QQ}$ is ample on $B_{\QQ}$, 
and $\adeg \left( \acherncl_1(\overline{L})^{d+1} \right) > 0$,
then $\overline{L}$ is big.
\end{Theorem}

Moreover, we have the following.

\begin{Proposition}
\label{prop:equiv:cond:big}
Let $B$ be a projective arithmetic variety with $d = \dim B_{\QQ}$, and
$\overline{L}$ a $C^{\infty}$-hermitian line bundle on $B$.
Then, the following are equivalent.
\begin{enumerate}
\renewcommand{\labelenumi}{(\arabic{enumi})}
\item
$\overline{L}$ is big.

\item
For any $C^{\infty}$-hermitian line bundles $\overline{M}$ on $B$,
there are a positive integer $n$ and
a non-zero section $s$ of $H^0(B, L^{\otimes n} \otimes M)$ with
$\Vert s \Vert_{\sup} < 1$.
\end{enumerate}
\end{Proposition}

\Proof
First, we assume (1).
Then, there is a non-zero section $s_1$ of $H^0(B, L^{\otimes n_1})$ with
$\Vert s_1 \Vert_{\sup} < 1$ for some $n_1$.
Moreover, since $\rank H^0(B, L^{\otimes m}) = O(m^d)$,
we can find a non-zero section $s_2$ of $H^0(B, L^{\otimes n_2} \otimes M)$.
Let $n_3$ be a sufficiently large integer with
\[
 \left(\Vert s_1 \Vert_{\sup}\right)^{n_3} \Vert s_2 \Vert_{\sup} < 1.
\]
Then, $s_1^{\otimes n_3} \otimes s_2 \in H^0(B, L^{\otimes n_3 n_1 + n_2} \otimes M)$, and
\[
\Vert s_1^{\otimes n_3} \otimes s_2 \Vert_{\sup} \leq 
\left(\Vert s_1 \Vert_{\sup}\right)^{n_3} \Vert s_2 \Vert_{\sup} < 1.
\]
Thus, we get (2).

\medskip
Next, we assume (2).
It is sufficient to show that $\rank H^0(B, L^{\otimes m}) = O(m^d)$.
Let $\overline{A}$ be an ample $C^{\infty}$-hermitian line bundle on $B$.
Then, there are a positive integer $n_1$ and a non-zero section
$s$ of $H^0(B, L^{\otimes n_1} \otimes A^{\otimes -1})$.
Thus, we have an injection $A \to L^{\otimes n_1}$. 
Therefore, $\rank H^0(B, L^{\otimes m}) = O(m^d)$.
\QED

Let $\aPic_{C^{\infty}}(B)$ (resp. $\aPic_{C^0}(B)$)be 
the isomorphic classes of
$C^{\infty}$-hermitian (resp. continuous hermitian)
line bundles on $B$.
An element of $\aPic_{C^{\infty}}(B) \otimes \QQ$
(resp. $\aPic_{C^0}(B) \otimes \QQ$) is called
a $C^{\infty}$-hermitian (resp. continuous hermitian)
$\QQ$-line bundle on $B$.
For simplicity, the group structure of $\aPic_{C^{\infty}}(B) \otimes \QQ$
(or $\aPic_{C^0}(B) \otimes \QQ$) is often written additively.
Note that the previous definitions of
`ample', `vertically nef', `horizontally nef', `nef' and `big'
work for $C^{\infty}$-hermitian $\QQ$-line bundles on $B$.
Let $\overline{L}$ be a continuous hermitian $\QQ$-line bundle.
We say $\overline{L}$ is {\em effective} if 
$\overline{L} \in \Pic_{C^0}(B)$ and
$H^0(B, L)$ contains a non-zero section $s$
with $\Vert s \Vert_{\sup} \leq 1$.
Moreover, $\overline{L}$ is said to be {\em $\QQ$-effective}
if $n \overline{L}$ is effective for some positive integer $n$.

\begin{Proposition}
\label{prop:intersection:nef:line:bundle}
Let $B$ be a projective arithmetic variety with $d = \dim B_{\QQ}$.
Then, we have the following.
\begin{enumerate}
\renewcommand{\labelenumi}{(\arabic{enumi})}
\item
If $\overline{L}_1, \ldots, \overline{L}_{d+1}$ are nef
$C^{\infty}$-hermitian $\QQ$-line bundles, then
\[
\adeg \left( \acherncl_1(\overline{L}_1) \cdots 
\acherncl_1(\overline{L}_{d+1})
\right) \geq 0.
\]

\item
If $\overline{L}_1, \ldots, \overline{L}_{d}$ are nef
$C^{\infty}$-hermitian $\QQ$-line bundles and $\overline{M}$ is a
$\QQ$-effective continuous hermitian $\QQ$-line bundle, then
\[
\adeg \left( \acherncl_1(\overline{L}_1) \cdots 
\acherncl_1(\overline{L}_{d}) \cdot \acherncl_1(\overline{M})
\right) \geq 0.
\]
\end{enumerate}
\end{Proposition}

\Proof
(1)
It can be proved by using Nakai-Moishezon's criterion on arithmetic varieties
(cf. \cite[Corollary~4.8]{ZhPos}).
Here we would like to give a more elementary proof,
which is a simpler case of Theorem~\ref{thm:estimate:liminf:inf:height}.
Let us begin with the following lemma.

\begin{Lemma}
\label{lem:intersection:ample:nef}
Let $\pi : B \to \Spec(\ZZ)$ be a projective arithmetic variety
with $d = \dim (B_{\QQ})$, and
$\overline{L}$ a nef $C^{\infty}$-hermitian $\QQ$-line bundle on $B$.
Moreover, let $\overline{A}$ be 
a vertically nef $C^{\infty}$-hermitian $\QQ$-line bundle such that
$A_{\QQ}$ is ample on $B_{\QQ}$ and,
for all integral subschemes $\Gamma$ of $B$ with $\pi(\Gamma) = \Spec(\ZZ)$,
\[
\adeg \left( \acherncl_1(\rest{\overline{A}}{\Gamma})^{\dim(\Gamma_{\QQ})+1} 
\right)
> 0.
\]
Then, for all $0 \leq i \leq d+1$,
\[
\adeg \left( \acherncl_1(\overline{L})^i \cdot
\acherncl_1(\overline{A})^{d+1-i} 
\right) \geq 0.
\]
\end{Lemma}

\Proof
We prove this lemma by induction on $d$.
If $d= 0$, then our assertion is obvious, so that
we assume that $d > 0$.

Case $i = 0, \ldots, d$:\quad
Since $\adeg \left( \acherncl_1(\overline{A})^{d+1} \right) > 0$,
replacing $\overline{A}$ by $\overline{A}^{\otimes n}$ ($n > 0$),
we may assume that there is a non-zero section $s$ of $H^0(B, A)$ with
$\Vert s \Vert_{\sup} < 1$.
Let $\zero(s) = a_1 \Gamma_1 + \cdots + a_e \Gamma_e$
be the decomposition of $\zero(s)$ as cycles.
Here, $a_j > 0$. Then,
\begin{multline*}
\adeg \left( \acherncl_1(\overline{L})^i \cdot
\acherncl_1(\overline{A})^{d+1-i} 
\right) = \sum_j a_j \adeg \left( \acherncl_1(\rest{\overline{L}}{\Gamma_j})^i \cdot
\acherncl_1(\rest{\overline{A}}{\Gamma_j})^{d-i} 
\right) \\
- \int_{B(\CC)}\log (\Vert s \Vert) c_1(\overline{L})^{\wedge i} \wedge
\acherncl_1(\overline{A})^{\wedge d-i}.
\end{multline*}
Therefore, by Lemma~\ref{lem:non:negative:prod:semi:positive}
and the hypothesis of induction, the above is non-negative.

Case $i = d+1$:\quad
Let $P(t)$ be a polynomial given by
\[
P(t) = \adeg \left( \acherncl_1( \overline{L} + t \overline{A})^{d+1} \right).
\]
Here, we claim the following.

\begin{Claim}
\label{claim:lem:intersection:ample:nef}
If $t > 0$and $P(t) > 0$, then
$P(t) \geq t^{d+1} \adeg \left( \acherncl_1(\overline{A})^{d+1} \right)$.
\end{Claim}

By using the hypothesis of induction and the assumption $P(t) > 0$,
we can see 
\[
\adeg \left( \acherncl_1(\rest{ (\overline{L} + t \overline{A})}{\Gamma})^{\dim(\Gamma_{\QQ})+1} \right)
> 0
\]
for all integral subschemes $\Gamma$ on $B$ with $\pi(\Gamma) = \Spec(\ZZ)$.
Thus, in the same way as above, we have
$\adeg \left( \acherncl_1(\overline{L}) \cdot
\acherncl_1(\overline{L} + t \overline{A})^{d} 
\right) \geq 0$. Hence,
\begin{align*}
P(t) & = \adeg \left( (\acherncl_1(\overline{L}) + t
\acherncl_1(\overline{A})) \cdot \acherncl_1(\overline{L} + t \overline{A})^d \right) \\
& \geq
t \adeg \left( \acherncl_1(\overline{A}) \cdot 
\acherncl_1(\overline{L} + t \overline{A})^d \right) \\
& \geq t^{d+1} \adeg \left( \acherncl_1(\overline{A})^{d+1} \right)
\end{align*}
Therefore, we get the claim.

\medskip
We set 
$t_0 = \max \{ t \in \RR \mid P(t) = 0 \}$.
Here we assume that $t_0 > 0$.
Then, by the above claim,
for all $t > t_0$,
\[
P(t) \geq t^{d+1} \adeg \left(
\acherncl_1(\overline{A})^{d+1} \right).
\]
Hence, taking $t \to t_0$,
\[
0 = P(t_0) \geq t_0^{d+1} \adeg \left(
\acherncl_1(\overline{A})^{d+1} \right) > 0.
\]
This is a contradiction, namely, $t_0 \leq 0$.
Thus, $P(0) = \adeg \left( \acherncl_1(\overline{L})^{d+1} \right) \geq 0$.
\QED

\bigskip
Let us go back to the proof of (1) of
Proposition~\ref{prop:intersection:nef:line:bundle}.
We prove this by induction on $d$.
Let $\overline{A}$ be an ample $C^{\infty}$-hermitian line bundle on $B$.
Let $\epsilon$ is a positive rational number.
Then, by Lemma~\ref{lem:intersection:ample:nef},
we can see
\[
\adeg \left(  \acherncl_1(\overline{L}_{d+1} + \epsilon \overline{A})^{d+1} \right)
> 0.
\]
Hence, using a small section of a positive multiple of
$\overline{L}_{d+1} + \epsilon \overline{A}$, the hypothesis of induction and
Lemma~\ref{lem:non:negative:prod:semi:positive},
we can see 
\[
\adeg \left( \acherncl_1(\overline{L}_1) \cdots \acherncl_1(\overline{L}_{d})
\cdot ( \overline{L}_{d+1} + \epsilon \overline{A} )
\right) \geq 0.
\]
Thus, we have our assertion taking $\epsilon \to 0$.

\medskip
(2)
This is a consequence of (1).
\QED

Finally let us consider the following lemma, which was used
in the proof of the previous proposition.

\begin{Lemma}
\label{lem:non:negative:prod:semi:positive}
Let $V$ be a vector space over $\RR$ with a complex structure $J$, i.e.,
an endomorphism $J : V \to V$ with $J^2 = -\operatorname{id}_V$.
Let $T$ \rom{(}resp. $T'$\rom{)} be the eigenspace of $J$ with respect to $\sqrt{-1}$
\rom{(}resp. $-\sqrt{-1}$\rom{)} in $V \otimes_{\RR} \CC$.
Note that the complex conjugation in $V \otimes_{\RR} \CC$ gives rise
to the anti-isomorphism of $T$ and $T'$ over $\CC$.
Let us fix a basis $\{ e_1, \ldots, e_n \}$ of $T$ over $\CC$.
For a hermitian $n \times n$-matrix $H = (h_{ij})$, we set
\[
\omega(H) = \sqrt{-1} \sum_{i,j} h_{ij} (e_i \wedge \bar{e}_j).
\]
If $H_1, \ldots, H_n$ are semipositive hermitian $n \times n$-matrices,
then there is a non-negative real number $\lambda$ with
\[
\omega(H_1) \wedge \cdots \wedge \omega(H_n) = \lambda
(\sqrt{-1})^n 
(e_1 \wedge \bar{e}_1) \wedge \cdots \wedge (e_n \wedge \bar{e}_n).
\]
\end{Lemma}

\Proof
First we claim the following.

\begin{Claim}
\label{claim:lem:non:negative:prod:semi:positive}
Let $x_1, \ldots, x_n$ be elements of $T$. If we set
$x_i = \sum_j a_{ij}e_j$ and $A = (a_{ij})$, then
\[
(x_1 \wedge \bar{x}_1) \wedge \cdots \wedge (x_n \wedge \bar{x}_n) =
\vert \det(A) \vert^2
(e_1 \wedge \bar{e}_1) \wedge \cdots \wedge (e_n \wedge \bar{e}_n).
\]
\end{Claim}

Since
\[
x_1 \wedge \cdots \wedge x_n = \det(A)(e_1 \wedge \cdots \wedge e_n)
\quad\text{and}\quad
\bar{x}_1 \wedge \cdots \wedge \bar{x}_n = \overline{\det(A)} 
(\bar{e}_1 \wedge \cdots \wedge \bar{e}_n),
\]
we have
\[
(x_1 \wedge \cdots \wedge x_n) \wedge 
(\bar{x}_1 \wedge \cdots \wedge \bar{x}_n)
= \vert \det(A) \vert^2 
(e_1 \wedge \cdots \wedge e_n) \wedge 
(\bar{e}_1 \wedge \cdots \wedge \bar{e}_n),
\]
which shows us the claim.

\bigskip
By our assumption, there are unitary matrices $U_i$'s and
non-negative real numbers $\lambda^i_1, \ldots, \lambda^i_n$ ($i = 1, \ldots, n$) such that
$U_i^* H_i U_i = \operatorname{diag}(\lambda^i_1, \ldots, \lambda^i_n)$ for all $i$,
where $A^* = \bar{A}^t$. Thus, if we take a new basis
$e^i_1, \ldots, e^i_n$ according to $U_i$,
then
\[
\omega(H_i) = \sqrt{-1} \sum_{a=1}^n \lambda^i_a (e^i_a \wedge \bar{e}^i_a).
\]
Thus, we obtain
\[
\omega(H_1) \wedge \cdots \wedge \omega(H_n) =
(\sqrt{-1})^n \sum_{1 \leq a_1, \ldots, a_n \leq n}
\lambda^1_{a_1} \cdots \lambda^n_{a_n}
(e^1_{a_1} \wedge \bar{e}^1_{a_1}) \wedge \cdots \wedge (e^n_{a_n} \wedge \bar{e}^n_{a_n}).
\]
On the other hand, by the above claim, there is a non-negative real number
$\tau_{a_1, \ldots, a_n}$ with
\[
(e^1_{a_1} \wedge \bar{e}^1_{a_1}) \wedge \cdots \wedge (e^n_{a_n} \wedge \bar{e}^n_{a_n})
= \tau_{a_1, \ldots, a_n}
(e_1 \wedge \bar{e}_1) \wedge \cdots \wedge (e_n \wedge \bar{e}_n).
\]
Therefore,
\[
\omega(H_1) \wedge \cdots \wedge \omega(H_n) =
\left(
\sum_{1 \leq a_1, \ldots, a_n \leq n}
 \tau_{a_1, \ldots, a_n} \lambda^1_{a_1} \cdots \lambda^n_{a_n}
\right)
(\sqrt{-1})^n 
(e_1 \wedge \bar{e}_1) \wedge \cdots \wedge (e_n \wedge \bar{e}_n).
\]
Hence, we get our lemma.
\QED

\renewcommand{\theTheorem}{\arabic{section}.\arabic{subsection}.\arabic{Theorem}}
\renewcommand{\theClaim}{\arabic{section}.\arabic{subsection}.\arabic{Theorem}.\arabic{Claim}}
\renewcommand{\theequation}{\arabic{section}.\arabic{subsection}.\arabic{Theorem}.\arabic{Claim}}
\section{Arithmetic height functions over finitely generated fields}
\subsection{Polarization of finitely generated fields over $\QQ$}
\setcounter{Theorem}{0}
Let $K$ be a finitely generated field over $\QQ$
with $\trdeg_{\QQ}(K) = d$.
A normal projective arithmetic variety $B$ is called
an {\em arithmetic model of $K$} if
the function field of $B$ is isomorphic to $K$.
A collection $(B; \overline{H}_1, \ldots, \overline{H}_d)$
of the arithmetic model $B$ of $K$ and
nef $C^{\infty}$-hermitian $\QQ$-line bundles $\overline{H}_1, \ldots, \overline{H}_d$
on $B$ is called a {\em polarization of $K$}.
Note that if $d = 0$, then we do not require any kind of $C^{\infty}$-hermitian line bundles
to fix a polarization of $K$.
For short, the polarization $(B; \overline{H}_1, \ldots, \overline{H}_d)$ is often
denoted by $\overline{B}$.
Moreover, the polarization $\overline{B}$ is said to be
{\em big} if $\overline{H}_1, \ldots, \overline{H}_d$ are big.
If $\overline{H}_1 = \cdots = \overline{H}_d$, say $\overline{H}$,
then the polarization $\overline{B}$ is 
simply called a {\em polarization of $K$
given by $\overline{H}$}.

Let $K'$ be a finite extension field of $K$, and
$\mu : B' \to B$ the normalization of $B$ in $K'$.
Then, we have a polarization 
$(B'; \mu^*(\overline{H}_1), \ldots, \mu^*(\overline{H}_d))$ of $K'$.
This polarization is denoted by $\overline{B}_{K'}$, and
is called the {\em polarization of $K'$ induced by $\overline{B}$}.
Clearly, if $\overline{B}$ is big, then
so is $\overline{B}_{K'}$.

Here let us consider the existence of a special polarization.

\begin{Proposition}
\label{prop:exist:line:pencil:type}
Let $K$ be a finitely generated field over $\QQ$ with $\trdeg_{\QQ} K = d$.
Then, there are a finite extension field $K'$ of $K$,
an arithmetic model $B$ of $K'$, and
a nef $C^{\infty}$-hermitian line bundle $\overline{H}$ on $B$ 
such that $H$ is ample
and $\adeg ( \acherncl_1(\overline{H})^{d+1} ) = 0$.
\end{Proposition}

\Proof
If $d = 0$, then we can take $\overline{H}$ as $(O_K, \vert\cdot\vert_{can})$,
where $O_K$ is the ring of integers in $K$.
Thus, we may assume that $d > 0$.

We first need a special arithmetic surface.
Let us consider the following elliptic curve due to J. Tate
(cf. \cite[5.10]{SerPG}):
\[
y^2 + xy + \epsilon^2 y = x^3,
\]
where $\epsilon = (5 + \sqrt{29})/2$ is
the fundamental unit of $\QQ(\sqrt{29})$. 
Then, the discriminant of this curve
is $-\epsilon^{10}$. 
We denote $\QQ(\sqrt{29})$ by $k$, and the ring of integers by $O_{k}$,
i.e., $O_k = \ZZ[\epsilon]$.
Here we set
\[
E = \Proj \left( O_{k}[X, Y, Z]/(Y^2Z + XYZ + \epsilon^2 YZ^2 - X^3) \right)
\]
and $E^d = E \times_{O_k} \cdots \times_{O_k} E$.
Then, since $E$ is smooth over $O_{k}$,
$E^d$ is an abelian scheme over $O_k$.
For an ample line bundle $L$ on $E^d$, we set
$H_0 = [-1]^*(L) \otimes L$. Then,
$H_0$ is symmetric on each fiber of $E^d \to \Spec(O_{k})$.
Thus, $[2]^*(H_0) = H_0^{\otimes 4}$ because the class group of $k$ is trivial.
Moreover, on each infinite fiber, we give the cubic metric of $H_0$
with $[2]^*(\overline{H}_0) = \overline{H}_0^{\otimes 4}$
(cf. \cite{MoBMet}),
so that $c_1(\overline{H}_0)$ is positive on each infinite fiber.
Thus, the height function $h_{\overline{H}_0}$ given by $\overline{H}_0$ is
nothing more than the N\'{e}ron-Tate height associated with $(H_0)_{k}$.
Hence, we can see that $\overline{H}_0$ is nef and
$\adeg ( \acherncl_1(\overline{H}_0)^{d+1} ) = 0$
by virtue of \cite[Theorem~(5.2)]{ZhPos}.

Let $K_0$ be the function field of $E^d$. Then,
$(E^d, \overline{H}_0)$ is a polarization of $K_0$.
Here we take a finite extension $K'$ of $K$ with $K_0 \subseteq K'$.
Then, the polarization of $K'$ induced by $(E^d, \overline{H}_0)$
is our desired polarization.
\QED

\subsection{Naive height functions}
\label{subsec:naive:ht:fun}
\setcounter{Theorem}{0}
Let $K$ be a finitely generated field over $\QQ$
with $\trdeg_{\QQ}(K) = d$, and
$\overline{B} = (B; \overline{H}_1, \ldots, \overline{H}_d)$ a polarization of $K$.
For $(x_0, \ldots, x_n) \in K^{n+1} \setminus \{ 0 \}$,
we set
\begin{multline*}
h^{\overline{B}}_{nv, K}(x_0, \ldots, x_n) =
\sum_{\Gamma}
\max_{i} \{ - \ord_{\Gamma}(x_i) \} \adeg 
\left( \acherncl_1(\rest{\overline{H}_1}{\Gamma}) \cdots 
\acherncl_1(\rest{\overline{H}_d}{\Gamma}) \right) \\
+ \int_{B(\CC)} \log \left( \max_i \{ |x_i| \} \right) 
c_1(\overline{H}_1) \wedge \cdots \wedge c_1(\overline{H}_d),
\end{multline*}
where $\Gamma$ runs over all prime divisors on $B$.
Note that if $d=0$, then the term
$ \acherncl_1(\rest{\overline{H}_1}{\Gamma}) \cdots 
\acherncl_1(\rest{\overline{H}_d}{\Gamma})$ is $[\Gamma]$ as cycle, and
$c_1(\overline{H}_1) \wedge \cdots \wedge c_1(\overline{H}_d)$ is $1$,
so that in this case, the above naive height coincides with the usual naive height
over a number field.

For $x \in K \setminus \{ 0 \}$,
\begin{align*}
0 & = \adeg \left( \acherncl_1(\overline{H}_1) \cdots \acherncl_1(\overline{H}_d)
\cdot \widehat{(x^{-1})} \right) \\
& =
\sum_{\Gamma} (- \ord_{\Gamma}(x)) \adeg 
\left( \acherncl_1(\rest{\overline{H}_1}{\Gamma}) \cdots 
\acherncl_1(\rest{\overline{H}_d}{\Gamma}) \right)
+ 
\int_{B(\CC)} \log (|x|) c_1(\overline{H}_1) \wedge \cdots \wedge c_1(\overline{H}_d).
\end{align*}
Thus, we can see that
\[
h^{\overline{B}}_{nv,K}(ax_0, \ldots, ax_n) = 
h^{\overline{B}}_{nv,K}(x_0, \ldots, x_n)
\]
for all $(x_0, \ldots, x_n) \in K^{n+1} \setminus \{ 0 \}$ and
all $a \in K \setminus \{ 0 \}$.
Hence, we have a function
\[
h^{\overline{B}}_{nv,K} : \PP^n(K) \to \RR.
\]

Let $K'$ be a finite extension field of $K$, and
$\overline{B}'$ the polarization of $K'$ induced by
$\overline{B}$.
Then, for $(x_0, \ldots, x_n) \in K^{n+1} \setminus \{ 0 \}$,
it is easy to see that
\[
h^{\overline{B}'}_{nv,K'}(x_0, \ldots, x_n) = 
[K' : K]h^{\overline{B}}_{nv,K}(x_0, \ldots, x_n).
\]
(Of course, this can be checked directly. In the next subsection, we give
an alternative definition of
$h^{\overline{B}}_{nv,K}$ in terms of Arakelov intersection theory, which
also shows the above formula by virtue of the projection formula.)
Thus, a family
${\displaystyle \left\{ [K' : K]^{-1} h^{\overline{B}'}_{nv,K'} 
\right\}_{K'}}$ of functions
gives rise to a naive height function 
\[
h^{\overline{B}}_{nv} : \PP^n(\overline{K}) \to \RR
\]
associated with $\overline{B}$.

\subsection{Height functions in terms of Arakelov intersection theory}
\setcounter{Theorem}{0}
Let $K$ be a finitely generated field over $\QQ$ 
with $\trdeg_{\QQ}(K) = d$, and
$\overline{B} = (B; \overline{H}_1, \ldots, \overline{H}_d)$ a polarization of $K$.
Let $X$ be a projective variety over $K$, and
$L$ a line bundle on $X$.
Let $f : \XX \to B$ be a projective morphism of arithmetic varieties, and
$\overline{\LL}$ a continuous hermitian $\QQ$-line bundle on $\XX$ such that
$\XX_K = X$ and $\LL_K = L$.
We say a pair $(\XX, \overline{\LL})$ is called a 
{\em model of $(X, L)$} over $B$.
Moreover, if $\overline{\LL}$ is a $C^{\infty}$-hermitian $\QQ$-line
bundle, then the pair $(\XX, \overline{\LL})$ is called a
{\em $C^{\infty}$-model} of $(X, L)$.
For $P \in X(\overline{K})$, we denote by $\Delta_P$ the Zariski
closure of $\Image(\Spec(\overline{K}) \overset{P}{\longrightarrow} \XX)$ in $\XX$.
Then, we define the height of $P$ with respect to
$(\XX, \overline{\LL})$ and $\overline{B}$ to be
\[
h_{(\XX, \overline{\LL})}^{\overline{B}}(P) =
\frac{1}{[K(P) : K]}
\adeg \left( \acherncl_1(\rest{f^* \overline{H}_1}{\Delta_P}) \cdots
\acherncl_1(\rest{f^* \overline{H}_d}{\Delta_P}) \cdot 
\acherncl_1(\rest{\overline{\LL}}{\Delta_P}) \right),
\]
where if $d=0$, then the term
$\acherncl_1(\rest{f^* \overline{H}_1}{\Delta_P}) \cdots
\acherncl_1(\rest{f^* \overline{H}_d}{\Delta_P})$
should be $[\Delta_P]$ as cycles.
\begin{comment}
In other word,
\begin{align*}
h_{(\XX, \overline{\LL})}^{\overline{B}}(P) & =
\frac{1}{[K(P) : K]}
\adeg \left( \acherncl_1(f^* \overline{H})^{d} \cdot 
\acherncl_1(\rest{\overline{\LL}}{\Delta_P}) \right) \\
& =
\frac{1}{[K(P) : K]}
\adeg \left( \acherncl_1(\overline{H})^d \cdot 
f_* (\acherncl_1(\rest{\overline{\LL}}{\Delta_P})) \right)
\end{align*}
\end{comment}
Let us begin with the following proposition.

\begin{Proposition}
\label{prop:base:change:field}
Let $K'$ be a finite extension field of $K$, and
let $\pi : B' \to B$ be a morphism of projective normal arithmetic
varieties such that the function field of $B'$ is $K'$.
Let $\XX'$ be the main component of $\XX \times_B B'$.
We set the induced morphism as follows:
\[
\begin{CD}
\XX @<{\pi'}<< \XX' \\
@V{f}VV @VV{f'}V \\
B @<{\pi}<< B'.
\end{CD} 
\]
Then, $h_{(\XX', {\pi'}^*(\overline{\LL}))}^{(B'; \pi^*(\overline{H}_1), \ldots,
\pi^*(\overline{H})_d))} =
[K': K]h_{(\XX, \overline{\LL})}^{(B; \overline{H}_1, \ldots, \overline{H}_d)}$.
\end{Proposition}

\Proof
Pick up $P \in X(\overline{K})$.
Let $\Delta_P$ (resp. $\Delta_P'$) be
the closure in $\XX$ (resp. $\XX'$).
Then, by the projection formula
(cf. (2) of Proposition~\ref{prop:projection:formula:intersection:num}),
\begin{align*}
h_{(\XX', {\pi'}^*(\overline{\LL}))}^{(B'; \pi^*(\overline{H}_1), \ldots,
\pi^*(\overline{H})_d))}(P) &
= \frac{\adeg \left( \acherncl_1(\rest{{\pi'}^* f^* \overline{H}_1}{\Delta_P'})
\cdots \acherncl_1(\rest{{\pi'}^* f^* \overline{H}_d}{\Delta_P'}) 
\cdot 
\acherncl_1(\rest{{\pi'}^* \overline{\LL}}{\Delta_P'}) \right)}{\deg(\Delta_P' \to B')} \\
& =
\frac{\deg(\Delta_P' \to \Delta_P)}{\deg(\Delta_P' \to B')}
\adeg \left( \acherncl_1(\rest{f^* \overline{H}_1}{\Delta_P}) \cdots
\acherncl_1(\rest{f^* \overline{H}_d}{\Delta_P}) \cdot 
\acherncl_1(\rest{\overline{\LL}}{\Delta_P}) \right) \\
& =
[K' : K]h_{(\XX, \overline{\LL})}^{(B; \overline{H}_1, \ldots, \overline{H}_d)}(P).
\end{align*}
\QED

Let $\PP^n_{\CC}$ be the $n$-dimensional projective space over $\CC$, and
$\OO(1)$ the tautological line bundle on $\PP^n_{\CC}$.
We fix a homogeneous coordinate $(X_0, \ldots, X_n)$ of $\PP^n_{\CC}$,
i.e., a basis of $H^0(\PP^n_{\CC}, \OO(1))$.
For a real number $l \geq 1$, we define the metric $\normabb_l$
of $\OO(1)$ to be
\[
\Vert X_i \Vert_{l} = \frac{|X_i|}{\left(|X_0|^l + \cdots + |X_n|^l\right)^{1/l}}.
\]
$\normabb_{2}$ is called the {\em Fubini-Study metric} and
is denoted by $\normabb_{FS}$.
Moreover, the metric $\normabb_{\infty}$ is defined by
\[
\Vert X_i \Vert_{\infty} = \frac{|X_i|}{\max \{|X_0|, \ldots, |X_n| \}}. 
\]
Note that $\lim_{l\to\infty}\normabb_l = \normabb_{\infty}$.

Let us consider $\PP^n_B$ and the natural projection
$p : \PP^n_B \to \PP^n_{\ZZ}$.
By abuse of notation, $p^*(\OO(1), \normabb_{l})$ is
denoted by $(\OO(1), \normabb_{l})$ for $1 \leq l \leq \infty$.
Then, we have the following proposition.

\begin{Proposition}
\label{prop:naive:height:intersection:height}
Let $h^{\overline{B}}_{nv}$ be the naive height on $\PP^n(\overline{K})$
defined in the
previous subsection~\rom{\ref{subsec:naive:ht:fun}}. Then, 
$h^{\overline{B}}_{nv}$
coincides with $h^{\overline{B}}_{(\PP^n_B, (\OO(1), \Vert\cdot\Vert_{\infty}))}$.
\end{Proposition}

\Proof
By virtue of Proposition~\ref{prop:base:change:field}
(actually in the same way as the proof of it),
it is sufficient to show the following claim.

\begin{Claim}
\label{claim:prop:naive:height:intersection:height}$
h^{\overline{B}}_{(\PP^n_B, (\OO(1), \Vert\cdot\Vert_{\infty}))}(P)
= h^{\overline{B}}_{nv}(P)$ for all $P \in X(K)$.
\end{Claim}

We first fix a basis $\{ X_0, \ldots, X_n \}$ of
$H^0(\PP^n_{\ZZ}, \OO(1))$.
Let $\Delta_P$ be the section corresponding with $P$, and
$s_P : B \to \Delta_P \hookrightarrow \PP^n_B$ the canonical morphism.
For simplicity, we assume that $s_P^*(X_0) \not= 0$.
If we set $a_i = s_P^*(X_i)/s_P^*(X_0)$ for $i=0,\ldots,n$,
then $a_i \in K$ and 
$h^{\overline{B}}_{nv}(P) = h^{\overline{B}}_{nv}(a_0, a_1, \ldots, a_n)$.
Since $\zero(s_P^*(X_0)) = \sum_{\Gamma}
\ord_{\Gamma}(s_P^*(X_0))\Gamma$,
\begin{multline*}
\adeg \left(
\acherncl_1(\overline{H})^d \cdot \acherncl_1(s_P^*(\OO(1),\normabb_{\infty}))
\right) =
\sum_{\Gamma}\ord_{\Gamma}(s_P^*(X_0))
\adeg\left( \acherncl_1(\rest{\overline{H}_1}{\Gamma}) \cdots \acherncl_1(\rest{\overline{H}_d}{\Gamma})
\right) \\
+ \int_{B(\CC)} 
\left( -\log s_P^*(\Vert X_0 \Vert_{\infty}) \right) c_1(\overline{H}_1) \wedge \cdots
\wedge c_1(\overline{H}_d)
\end{multline*}
By the definition of $\normabb_{\infty}$,
we can see that
$-\log s_P^*(\Vert X_0 \Vert_{\infty}) =  
\log \left( \max_i \{ |a_i| \} \right) $.
On the other hand,
since $s_P^*(X_i)$'s generate $s_P^*(\OO(1))$,
we have
\begin{align*}
\ord_{\Gamma}(s_P^*(X_0)) & =
\length_{\Gamma}\left( 
\frac{s_P^*(\OO(1))}{s_P^*(X_0)} \right) \\
& = \length_{\Gamma}\left( 
\frac{\OO_B s_P^*(X_0) + \cdots + \OO_B s_P^*(X_n)}
{s_P^*(X_0)} \right) \\
& =
\length_{\Gamma}\left( 
\frac{\OO_B a_0 + \cdots + \OO_B a_n}{\OO_B} 
\right) \\
& = \max_i \{ - \ord_{\Gamma}(a_i)\}.
\end{align*}
Thus, we get our claim.

Note that combining the above claim with 
Proposition~\ref{prop:base:change:field},
we can see that $h^{\overline{B}}_{nv,K} = [K' : K]h^{\overline{B}}_{nv,K'}$
as we remarked in the previous subsection~\ref{subsec:naive:ht:fun}.
\QED

Next let us consider the following proposition.

\begin{Proposition}
\label{prop:bound:below:base:locus}
If we denote $\Supp\left(\Coker(H^0(X, L) \otimes \OO_X \to L)\right)$
by $\Bs(L)$,
then there is a constant $C$ with
$h^{\overline{B}}_{(\XX, \overline{\LL})}(P) \geq C$
for all $P \in (X \setminus \Bs(L))(\overline{K})$.
\end{Proposition}

\Proof
Let $A$ be an ample line bundle on $B$ such that
$f_*(\LL) \otimes A$ is generated by global sections, i.e.,
\[
 H^0(B, f_*(\LL) \otimes A) \otimes \OO_B \to f_*(\LL) \otimes A
\]
is surjective. Hence,
$H^0(\XX, \LL \otimes f^*(A)) \otimes \OO_{\XX} \to 
f^*(f_*(\LL)) \otimes f^*(A)$ is surjective.
Therefore, if we set
\[
 \mathcal{Z} = \Supp\left(
 \Coker(H^0(\XX, \LL \otimes f^*(A)) \otimes \OO_{\XX} \to
        \LL \otimes f^*(A))\right),
\]
then $\mathcal{Z}_K = \Bs(L)$.

Let $\{ s_1, \ldots, s_l \}$ be a free basis of
$H^0(\XX, \LL \otimes f^*(A))$ as $\ZZ$-modules.
Let us choose a metric of $A$ such that
$\Vert s_i \Vert_{\sup} < 1$ for all $i$.
Pick up an arbitrary $P \in (X \setminus \Bs(L))(\overline{K})$.
Then, there is an $s_i$ with $s_i(P) \not= 0$,
which shows us that
$\rest{\overline{\LL} \otimes f^* \overline{A}}{\Delta_P}$
is effective. Thus, by (2) of Proposition~\ref{prop:intersection:nef:line:bundle},
\[
\adeg \left( \acherncl_1(\rest{f^* \overline{H}_1}{\Delta_P}) \cdots
\acherncl_1(\rest{f^* \overline{H}_d}{\Delta_P}) \cdot 
\acherncl_1(\rest{\overline{\LL} \otimes f^* \overline{A}}{\Delta_P}) 
\right)
\geq 0.
\]
Hence, by virtue of the projection formula
((2) of Proposition~\ref{prop:projection:formula:intersection:num}), 
we have
\begin{multline*}
\adeg \left( \acherncl_1(\rest{f^* \overline{H}_1}{\Delta_P}) \cdots
\acherncl_1(\rest{f^* \overline{H}_d}{\Delta_P}) \cdot 
\acherncl_1(\rest{\overline{\LL}}{\Delta_P}) \right) \\
+ \deg(\Delta_P \to B) \adeg(\acherncl_1(\overline{H}_1) \cdots
\acherncl_1(\overline{H}_d) \cdot 
\acherncl_1(\overline{A}))
\geq 0.
\end{multline*}
Therefore,
\[
h_{(\XX, \overline{\LL})}^{\overline{B}}(P) \geq 
- \adeg(\acherncl_1(\overline{H}_1) \cdots \acherncl_1(\overline{H}_d) \cdot \acherncl_1(\overline{A})).
\]
\QED

\begin{Corollary}
\label{cor:bound:trivial}
If $\LL_K = \OO_X$, then there is a constant $C$
with $\vert h^{\overline{B}}_{(\XX, \overline{\LL})}(P) \vert \leq C$
for all $P \in X(\overline{K})$.
\end{Corollary}

\Proof
Apply Proposition~\ref{prop:bound:below:base:locus} to
$\overline{\LL}$ and $\overline{\LL}^{\otimes -1}$.
\QED

\begin{Corollary}
\label{cor:height:well:def:const}
Let $(\XX, \overline{\LL})$ and $(\XX', \overline{\LL}')$ be two
models of $(X, L)$. Then, there is a constant $C > 0$ with
\[
\vert h^{\overline{B}}_{(\XX, \overline{\LL})}(P) - 
h^{\overline{B}}_{(\XX', \overline{\LL}')}(P)
\vert \leq C
\]
for all $P \in X(\overline{K})$.
\end{Corollary}

\Proof
Let us consider the graph $\XX''$ of the birational map $\XX \dasharrow \XX'$.
Let $\mu : \XX'' \to \XX$ and $\mu' : \XX'' \to \XX'$ 
be the canonical morphisms.
Then, by the projection formula, we can see that
\[
h^{\overline{B}}_{(\XX'', \mu^*(\overline{\LL}))} = 
h^{\overline{B}}_{(\XX, \overline{\LL})}
\quad\text{and}\quad
h^{\overline{B}}_{(\XX'', (\mu')^*(\overline{\LL}'))} = 
h^{\overline{B}}_{(\XX', \overline{\LL}')}.
\]
On the other hand, since $\mu^*(\overline{\LL})$ coincides with
$(\mu')^*(\overline{\LL}')$ on the generic fiber of $\XX'' \to B$.
Thus, by Corollary~\ref{cor:bound:trivial},
there is a constant $C$ with
\[
\vert h^{\overline{B}}_{(\XX'', \mu^*(\overline{\LL}))} 
- h^{\overline{B}}_{(\XX'', (\mu')^*(\overline{\LL}'))} \vert \leq C.
\]
Therefore, we get our corollary.
\QED

\begin{Definition}
By the above corollary, 
the class of $h^{\overline{B}}_{(\XX, \overline{\LL})}$
modulo
the set of all bounded functions on $X(\overline{K})$
does not depend on the choice of the model $(\XX, \overline{\LL})$
of $(X, L)$ over $B$.
This class is called 
the {\em height associated with $L$ and 
$\overline{B}$}, and
is denoted by $h^{\overline{B}}_L$.
By abuse of notation,
we often view $h^{\overline{B}}_L$ as a representative
of $h^{\overline{B}}_L$.
\end{Definition}

Here we list elementary properties of
height functions.

\begin{Proposition}
\label{prop:elem:prop:height:fun}
\begin{enumerate}
\renewcommand{\labelenumi}{(\arabic{enumi})}
\item
If $X \subseteq \PP_K^n$ and $L = \rest{\OO_{\PP^n_K}(1)}{X}$,
then $h^{\overline{B}}_L = \rest{h^{\overline{B}}_{nv}}{X(\overline{K})} 
+ O(1)$.
\item
For line bundles $L$ and $M$ on $X$,
$h^{\overline{B}}_{L \otimes M} = h^{\overline{B}}_{L} +
h^{\overline{B}}_{M} + O(1)$ and
$h^{\overline{B}}_{L^{\otimes -1}} = -h^{\overline{B}}_{L} + O(1)$.

\item
$h^{\overline{B}}_L$ is bounded below
on $(X \setminus \SBs(L))(\overline{K})$,
where $\SBs = \bigcap_{n > 0} \Bs(L^{\otimes n})$.
In particular, we have the following.
\begin{enumerate}
\renewcommand{\labelenumii}{(\arabic{enumi}.\arabic{enumii})}
\item
If $L$ is ample, then $h^{\overline{B}}_L$ is bounded below.

\item
If $L = \OO_X$, then $h^{\overline{B}}_L = O(1)$.
\end{enumerate}

\item
\rom{(}Northcott's theorem for our height functions\rom{)}
If $ \overline{B} = (B; \overline{H}_1, \ldots, \overline{H}_d)$ is big, i.e.,
$\overline{H}_i$'s are nef and big, then,
for any numbers $M$ and any positive integers $e$, the set
\[
\{ P \in X(\overline{K}) \mid h^{\overline{B}}_L(P) \leq M,
\quad [K(P) : K] \leq e \} 
\]
is finite.

\item
Let $\overline{H}_1', \ldots, \overline{H}_d'$ be 
nef $C^{\infty}$-hermitian line bundles on $B$ 
such that $\overline{H}'_i \otimes \overline{H}_i^{\otimes -1}$
is $\QQ$-effective for every $i$.
If $L$ is ample, then
$h^{(B; \overline{H}_1, \ldots, \overline{H}_d)}_L \leq 
h^{(B; \overline{H}'_1, \ldots, \overline{H}'_d)}_L + O(1)$.
\end{enumerate}
\end{Proposition}

\Proof
(1): This is derived from 
Proposition~\ref{prop:naive:height:intersection:height}
and Corollary~\ref{cor:height:well:def:const}.

\medskip
(2):
This follows from the formulae:
$h^{\overline{B}}_{(\XX, \overline{\LL} \otimes \overline{\mathcal{M}})}
= h^{\overline{B}}_{(\XX, \overline{\LL})} + h^{\overline{B}}_{(\XX, \overline{\mathcal{M}})}$
and
$h^{\overline{B}}_{(\XX, \overline{\LL}^{\otimes -1})}
= -h^{\overline{B}}_{(\XX, \overline{\LL})}$.

\medskip
(3):
Since there is a positive integer $n$ with
$\SBs(L) = \Bs(L^{\otimes n})$,
it is a consequence of
(2) and Proposition~\ref{prop:bound:below:base:locus}.

\medskip
(4):
This will be proved in \S\ref{sec:northcott:thm}
(cf. Theorem~\ref{thm:northcott:thm:fun:field}).

\medskip
(5):
Clearly, we may assume that $L$ is very ample.
Let $\phi_L : X \to \PP^N_K$ be the embedding by $L$.
Let $\XX$ be the closure of $X$ in $\PP^N_B$, and
$\OO_{\XX}(1)$ the restriction of $\OO(1)$ on $\PP^N_B$.
Then, $\OO_{\XX}(1)$ is $f$-ample, where $f$ is the canonical morphism
$\XX \to B$. Thus, there is an ample line bundle $Q$ on $B$ such that
$\OO_{\XX}(1) \otimes f^*(Q)$ is ample.
We set $\LL = \OO_{\XX}(1) \otimes f^*(Q)$.
Then, $\LL_K = L$. Moreover, we give a $C^{\infty}$-hermitian metric
to $\LL$ such that $\overline{\LL}$ is ample.

Let us pick up $P \in X(\overline{K})$ and 
let $\Delta_P$ be the closure of $P$ in $\XX$.
For simplicity, we set $A_i = \rest{\acherncl_1(f^* \overline{H}'_i)}{\Delta_P}$,
$B_i = \rest{\acherncl_1(f^* \overline{H}_i)}{\Delta_P}$, and
$C = \rest{\acherncl_1(\overline{\LL})}{\Delta_P}$,
Then,
\[
A_1 \cdots A_d \cdot C - B_1 \cdots B_d \cdot C = 
\sum_{i=0}^{d-1} A_1 \cdots A_i \cdot(A_{i+1} - B_{i+1}) \cdot B_{i+2} \cdots B_{d} \cdot C.
\]
On the other hand, since $A_i$, $B_i$, and $C$ are nef, and $A_{i+1} - B_{i+1}$
is $\QQ$-effective, by (2) of Proposition~\ref{prop:intersection:nef:line:bundle},
we have
\[
\adeg (A_{1} \cdots A_i \cdot(A_{i+1} - B_{i+1}) \cdot B_{i+2} \cdots B_d \cdot C) \geq 0.
\]
Thus, we get $\adeg(A_1 \cdots A_d \cdot C) \geq \adeg(B_1 \cdots B_d \cdot C)$, 
which says us that
\[
h^{(B; \overline{H}_1, \ldots, \overline{H}_d)}_{(\XX, \overline{\LL})}(P) \leq 
h^{(B; \overline{H}'_1, \ldots, \overline{H}'_d)}_{(\XX, \overline{\LL})}(P).
\]
Hence, we obtain (5).
\QED

\subsection{Canonical height functions on abelian varieties}
\setcounter{Theorem}{0}
Let $K$ be a finitely generated field over $\QQ$, 
and $\overline{B} = (B; \overline{H}_1, \ldots, \overline{H}_d)$
a polarization of $K$.
Let $A$ be an abelian variety over $K$, and
$L$ a line bundle on $A$.
Then, by virtue of the cubic theorem
and (3.2) of Proposition~\ref{prop:elem:prop:height:fun},
\[
h^{\overline{B}}_L(x + y + z)
-h^{\overline{B}}_L(x+y)
-h^{\overline{B}}_L(y+z)
-h^{\overline{B}}_L(z+x)
+h^{\overline{B}}_L(x)+h^{\overline{B}}_L(y)+h^{\overline{B}}_L(z)
\]
is a bounded function on $A(\overline{K}) \times A(\overline{K})
\times A(\overline{K})$.
Thus, there are a unique bilinear form
$q^{\overline{B}}_L : A(\overline{K}) \times A(\overline{K}) \to \RR$ and
a unique linear function $l^{\overline{B}}_L : A(\overline{K}) \to \RR$
such that
\[
h^{\overline{B}}_L(x) = q^{\overline{B}}_L(x,x) + l^{\overline{B}}_L(x) + 
O(1)
\]
(cf. \cite[Chapter~5, \S1]{LaFund}).
Actually, 
$q^{\overline{B}}_L(x,x)$ and $l^{\overline{B}}_L(x)$ are given by
the following formula:
\[
q^{\overline{B}}_L(x,x) = \lim_{n \to \infty} 2^{-2n}
h^{\overline{B}}_L(2^n x)
\quad\text{and}\quad
l^{\overline{B}}_L(x) = \lim_{n \to \infty} 2^{-n} \left(
2^{-2n} h^{\overline{B}}_L(2^n x) - q^{\overline{B}}_L(x,x) \right).
\]
$q^{\overline{B}}_L + l^{\overline{B}}_L$ is denoted by
$\hat{h}^{\overline{B}}_L$, and is called the {\em canonical height
function of $L$} with respect to the polarization $\overline{B}$.
Moreover, it is easy to see that $q = 0$ if $[-1]^*(L) = L^{\otimes -1}$, and
$l = 0$ if $[-1]^*(L) = L$.
Thus, if $L$ is symmetric, then
$\hat{h}^{\overline{B}}_L(x) = \lim_{n \to \infty} 2^{-2n}
h^{\overline{B}}_L(2^n x)$.
Here let us consider the following two propositions.

\begin{Proposition}
\label{prop:positivity:canonical:height}
Let $L$ be a symmetric ample line bundle on $A$.
Then, we have the following.
\begin{enumerate}
\renewcommand{\labelenumi}{(\arabic{enumi})}
\item
$\hat{h}^{\overline{B}}_L(x) \geq 0$ for all $x \in A(\overline{K})$.

\item
If $x$ is a torsion point, then $\hat{h}^{\overline{B}}_L(x) = 0$.

\item
We assume that $\overline{B}$ is big, i.e.,
$\overline{H}_1, \ldots, \overline{H}_d$ are nef and big.
Then,
$\hat{h}^{\overline{B}}_L(x) = 0$ if and only if
$x$ is a torsion point.
\end{enumerate}
\end{Proposition}

\Proof
(1)
This is a consequence of
(3.1) of Proposition~\ref{prop:elem:prop:height:fun}.

\medskip
(2)
We assume that $x$ is a torsion point.
Then, there is a positive number $n$ with $nx = 0$.
Thus, 
$0 = \hat{h}^{\overline{B}}_L(nx) = n^2 \hat{h}^{\overline{B}}_L(x)$.
Hence $\hat{h}^{\overline{B}}_L(x) = 0$.

\medskip
(3)
We assume that $\hat{h}^{\overline{B}}_L(x) = 0$.
Let us consider the subgroup $G$ generated by $x$.
If $x$ is defined over a finite extension field $K'$,
then every element of $G$ is defined over $K'$.
Moreover, the height of every element of $G$ is zero.
Thus, by (4) of Proposition~\ref{prop:elem:prop:height:fun},
$G$ is a finite group, namely,
$x$ is a torsion point.
\QED

\begin{Proposition}
\label{prop:comp:canonical:height}
Let $L$ and $M$ be symmetric line bundles on $A$.
Then we have the following.
\begin{enumerate}
\renewcommand{\labelenumi}{(\arabic{enumi})}
\item
$\hat{h}^{\overline{B}}_{L \otimes M} = \hat{h}^{\overline{B}}_L +
\hat{h}^{\overline{B}}_M$ and
$\hat{h}^{\overline{B}}_{L^{\otimes -1}} = -\hat{h}^{\overline{B}}_{L}$.

\item
If $L$ is ample, then there is a positive number $a$
with $\hat{h}^{\overline{B}}_M \leq a \hat{h}^{\overline{B}}_L$.

\item
Let $\overline{B}' = (B; \overline{H}'_1, \ldots, \overline{H}'_d)$ be 
another polarization of $K$.
If $L$ is ample and $\overline{B}'$ is big, then
there is a positive number $b$ with $\hat{h}^{\overline{B}}_L
\leq b \hat{h}^{\overline{B}'}_L$.
\end{enumerate}
\end{Proposition}

\Proof
(1) This is a consequence of (2) of
Proposition~\ref{prop:elem:prop:height:fun}.

\medskip
(2) There is a positive integer $a$ such that
$L^{\otimes a} \otimes M^{\otimes -1}$ is ample.
Thus, by (1) of Proposition~\ref{prop:positivity:canonical:height},
$\hat{h}^{\overline{B}}_{L^{\otimes a} \otimes M^{\otimes -1}} \geq 0$.
Hence, our assertion follows from (1).

\medskip
(3) By Proposition~\ref{prop:equiv:cond:big},
there are positive integers $n_1, \ldots, n_d$ such that
${\overline{H}'_i}^{\otimes n_i} \otimes \overline{H}^{\otimes -1}$
is effective for every $i$.
Thus it follows from (5) of Proposition~\ref{prop:elem:prop:height:fun}.
\QED

\subsection{Remarks}
\setcounter{Theorem}{0}
\begin{Remark}
Note that in general, 
our height function over a finitely generated field $K$ is
not a height function in the sense of Lang's book
\cite{LaFund}. For, the map $v : K \to \RR_{+}$ given by
\[
v(x) = \exp\left( \int_{B(\CC)} \log(\vert x \vert) c_1(\overline{H}_1)
\wedge \cdots \wedge c_1(\overline{H}_d) \right)
\]
is not necessarily a valuation of $K$,
where $(B; \overline{H}_1, \ldots, \overline{H}_d)$ is a polarization of $K$ and $d = \trdeg_{\QQ}(K)$.
\end{Remark}

\begin{Remark}
After writing the first draft of this paper,
Prof. Silverman kindly informed me the work of Altman.
In \cite{Alt}, he gave the size function similar to
our height functions. On an abelian variety $A$ over a field $K$,
he proved that there is a quadratic function $A(K) \to \RR$
with $\operatorname{size}(x) \leq Q(x)$ for all $x \in A(K)$.
Compared with his method, our way gives rise to
the point of view of geometry, so that
it is easy to handle it in the functorial framework.
\end{Remark}

\begin{Remark}
Here, we would like to point out a similarity between
our height functions and the characteristic function in Nevanlinna theory.
Let us choose $f \in \QQ(z) \setminus \{ 0 \}$.
If $f$ has no pole at $0$, then the characteristic function $T_f$ is
given by
\[
T_f(r) =
\sum_{|x|< r} \max\{ -\ord_x(f), 0 \} \log ( r/|x| )+ 
\int_{0}^{2\pi} \log^+ \left(| f(r e^{\sqrt{-1}\theta}) |\right) 
\frac{d\theta}{2\pi},
\]
where $\log^+(x) = \log(\max\{ x, 1 \})$.

On the other hand, if we fix a polarization 
$(\PP^1_{\ZZ}, (\OO(1), \normabb_{FS}))$ of $\QQ(z)$,
then the naive height of $(f : 1) \in \PP^1(\QQ(z))$
is given by
\begin{align*}
h_{nv}(f : 1) & =
\sum_{\Gamma}
\max \{ - \ord_{\Gamma}(f), 0 \} \adeg 
\left( \rest{(\OO(1), \normabb_{FS})}{\Gamma} \right) +
\int_{\CC} \log^+ (|f|)
\frac{\sqrt{-1} dz \wedge d \bar{z}}{2 \pi (1 + |z|^2)^2} \\
& = \sum_{\Gamma}
\max \{ - \ord_{\Gamma}(f), 0 \} \adeg 
\left( \rest{(\OO(1), \normabb_{FS})}{\Gamma} \right) \\
& \qquad\qquad\qquad +
\int_{0}^{\infty} \left[ \frac{2r}{(1 + |r|^2)^2}  \int_0^{2\pi} 
\log^+ \left(|f(r e^{\sqrt{-1}\theta})|\right) \frac{d \theta}{2\pi}
\right] dr.
\end{align*}
\end{Remark}
\renewcommand{\theTheorem}{\arabic{section}.\arabic{Theorem}}
\renewcommand{\theClaim}{\arabic{section}.\arabic{Theorem}.\arabic{Claim}}
\renewcommand{\theequation}{\arabic{section}.\arabic{Theorem}.\arabic{Claim}}

\section{Northcott's theorem over finitely generated fields}
\label{sec:northcott:thm}
Let $\CC[z_1, \ldots, z_n]$ be the ring of $n$-variables polynomials 
over $\CC$.
For $f \in \CC[z_1, \ldots, z_n]$, we denote by $|f|$
the maximal of the absolute values of coefficients of $f$.
Moreover, we denote by $\deg_i(f)$ the degree of $f$ with respect to $z_i$.

Let us consider the following $(1, 1)$-form $\omega_i$ on $\CC^n$ for each $i$:
\[
\omega_i = \frac{\sqrt{-1} dz_i \wedge d \bar{z}_i}{2 \pi (1 + |z_i|^2)^2}.
\]
For $f \in \CC[z_1, \ldots, z_n]$, we define $v(f)$ to be
\[
v(f) = \exp \left( \int_{\CC^n} \log (|f(z_1, \ldots, z_n)|) 
\omega_1 \wedge \cdots \wedge \omega_n \right).
\]
Then we have the following lemma.

\begin{Lemma}
\label{lem:estimate:value:integral}
For any $f \in \CC[z_1, \ldots, z_n]$,
$|f| \leq 2^{\deg_1(f) + \cdots + \deg_n(f)} v(f)$.
In particular, for any numbers $M$ and any non-negative integers
$d_1, \ldots, d_n$, the set
\[
\{ f \in \ZZ[z_1, \ldots, z_n] \mid v(f) \leq M, \quad
\deg_i(f) \leq d_i \  (i=1, \ldots, n) \}
\]
is finite.
\end{Lemma}

\Proof
First, let us consider the case $n=1$.
By straightforward calculations together with Jensen's formula, we can
see that $v(z-\alpha) = \sqrt{1 + |\alpha|^2}$ for all $\alpha \in \CC$.
Thus, if we set $f(z) = c(z-\alpha_1) \cdots (z-\alpha_d)$, then
$v(f) = \vert c \vert \sqrt{1 + |\alpha_1|^2} \cdots \sqrt{1 + |\alpha_d|^2}$.
Therefore, we can easily see that $|f| \leq 2^d v(f)$.

\medskip
In general, we will prove this lemma by induction on $n$.
We set
\[
f = a_0(z_2, \ldots, z_{n})z_1^{d_1} + a_1(z_2, \ldots, z_{n})z_1^{d_1 - 1} +
\cdots + a_{d_1}(z_2, \ldots, z_{n}),
\]
where $d_1 = \deg_1(f)$.
Then, by the case $n=1$,
\[
\max_i \{ |a_i(c_2, \ldots, c_{n})| \} \leq 2^{\deg_1(f)} \exp \left(
\int_{\CC} \log(|f(z_1, c_2, \cdots, c_{n})|) \omega_1 \right).
\]
for all $(c_2, \ldots, c_{n}) \in \CC^{n-1}$.
Thus, by hypothesis of induction,
\begin{align*}
\log(v(f)) & =
\int_{\CC^{n-1}} \left( \int_{\CC} \log (|f(z_1, \ldots, z_n)|) 
\omega_1 \right) 
\omega_{2} \wedge \cdots \wedge \omega_{n} \\
& \geq -\log(2)\deg_1(f) + \max_{a_i \not= 0}
\int_{\CC^{n-1}} \log(|a_i(z_2, \ldots, z_n)|) 
\omega_2 \wedge \cdots \wedge \omega_{n} \\
& \geq -\log(2)\deg_1(f) -\log(2)(\deg_2(f) + \cdots + \deg_n(f)) +
\max_{a_i \not = 0} \log |a_i| \\
& = -\log(2)(\deg_1(f) + \cdots + \deg_n(f)) + \log(|f|).
\end{align*}
Therefore, we get our lemma.
\QED

Next let us consider the following lemma.

\begin{Lemma}
\label{lem:northcott:special:case}
Let $B = (\PP^1_{\ZZ})^d$
and $\overline{H} = p_1^*((\OO_{\PP}(1), \normabb_{FS}))
\otimes \cdots \otimes p_d^*((\OO_{\PP}(1), \normabb_{FS}))$,
where $p_i$ is the projection to the $i$-th factor.
Then, for any numbers $M$, the set
\[
\{ P \in \PP^n(\QQ(z_1, \cdots, z_d)) \mid h_{nv}^{\overline{B}}(P) \leq M \}
\]
is finite, where $\overline{B} = (B; \overline{H}, \ldots, \overline{H})$
is a polarization given by $\overline{H}$.
\end{Lemma}

\Proof
Let $\overline{L}_1, \ldots, \overline{L}_d$ be $C^{\infty}$-hermitian
line bundles on $\PP^1_{\ZZ}$. We set 
$\overline{L} = p_1^*(\overline{L}_1) \otimes \cdots \otimes 
p_d^*(\overline{L}_d)$.
Let $\Delta_{\infty}$ be the closure of $\infty \in \PP^1_{\QQ}$
in $\PP^1_{\ZZ}$. We set $\Delta^{(i)}_{\infty} = p_i^*(\Delta_{\infty})$.
First of all, we claim the following:

\begin{Claim}
\label{claim:lem:northcott:special:case}
\begin{multline*}
\adeg(\acherncl_1(\rest{\overline{L}}{\Delta^{(i)}_{\infty}})^d) =\\ 
d! \left( \prod_{\substack{1 \leq k \leq d \\ k \not=i}} \deg(L_k) \right)
\adeg( \rest{\overline{L}_i}{\Delta_{\infty}} ) +
\frac{d!}{2} \sum_{\substack{1 \leq j \leq d \\ j \not= i}} 
\left( \prod_{\substack{1 \leq k \leq d \\ k \not=i,j}} \deg(L_k) \right)
\adeg(\acherncl_1(\overline{L}_j)^2).
\end{multline*}
\end{Claim}

For simplicity, we assume $i=d$.
Since $\Delta_{\infty} \simeq \Spec(\ZZ)$,
there is an isomorphism
\[
 \phi : \underbrace{\PP^1_{\ZZ} \times_{\ZZ} \cdots \times_{\ZZ}
\PP^1_{\ZZ}}_{\text{$(d-1)$-times}}
\overset{\sim}{\longrightarrow} \Delta_{\infty}^{(d)}.
\]
Let $q_i : (\PP^1_{\ZZ})^{d-1} \to \PP^1_{\ZZ}$ 
be the projection to the $i$-th factor.
Then, $\phi^*(p_i^*(\overline{L}_i)) = q_i^*(\overline{L}_i)$ for
$i = 1, \ldots, d-1$. 
Moreover, if we set $c =
\exp(-\adeg(\rest{\overline{L}_d}{\Delta_{\infty}}))$,
then $\phi^*(p_d^*(\overline{L}_d)) = (\OO, c \Vert \cdot \Vert_{can})$.
We need to calculate $\adeg(\acherncl_1
(\phi^*(\rest{\overline{L}}{\Delta^{(d)}_{\infty}}))^{d})$.
Since
\[
\acherncl_1(\phi^*(\rest{\overline{L}}{\Delta^{(d)}_{\infty}}))=
q_1^*(\acherncl_1(\overline{L}_1)) + \cdots +
q_{d-1}^*(\acherncl_1(\overline{L}_{d-1})) + (0, -2 \log(c)),
\]
it is equal to
\[
 \sum_{a_1 + \cdots + a_d = d}
\frac{d!}{a_1! \cdots a_d!}
\adeg \left(
q_1^*(\acherncl_1(\overline{L}_1))^{a_1} \cdots
q_{d-1}^*(\acherncl_1(\overline{L}_{d-1}))^{a_{d-1}} \cdot
(0, -2 \log(c))^{a_d}
\right).
\]
Let us consider non-zero terms in the above equation.
Clearly, $a_d$ must be $0$ or $1$. If $a_d = 0$, then
one of $a_1, \ldots, a_{d-1}$ is $2$, and others are $1$.
If $a_d = 1$, then $a_1 = \cdots = a_{d-1} = 1$.
Thus, $\adeg(\acherncl_1
(\phi^*(\rest{\overline{L}}{\Delta^{(d)}_{\infty}}))^{d})$
is equal to
\[
\frac{d!}{2} \sum_{j=1}^{d-1} 
\left( \prod_{\substack{1 \leq k \leq d-1 \\ k \not= j}}
\deg(L_k) \right)
\adeg(\acherncl_1(\overline{L}_j)^2) -
d! \log(c)\deg(L_1)\cdots\deg(L_{d-1}).
\]
Therefore, we get our claim.

\medskip
Let us go back to the proof of our lemma.
We fix a number $c$ with $0 < c < 1$.
We set
\[
 \overline{A}_i = p_1^*((\OO_{\PP}(1), \normabb_{FS}))
\otimes \cdots \otimes
p_i^*(\OO_{\PP}, c\normabb_{can})
\otimes \cdots \otimes  p_d^*((\OO_{\PP}(1), \normabb_{FS})) 
\]
Then, by the above claim, if we set
\[
 e =
d! \left(
-\log(c) + \frac{d-1}{2}\adeg(\acherncl_1(\OO(1), \normabb_{FS})^2)
\right),
\]
then
\[
\adeg \left( \acherncl_1(\rest{\overline{A}_i}{\Delta^{(j)}_{\infty}})
\right) =
\begin{cases}
e & \text{if $i=j$}, \\
0 & \text{if $i \not= j$}.
\end{cases} 
\]
Moreover, since $\overline{H}$ is ample,
there is a positive integer $n_0$ such that
$\overline{H}^{\otimes n_0} \otimes \overline{A}_i^{\otimes -1}$ is
effective for every $i$.
Thus, by (5) of Proposition~\ref{prop:elem:prop:height:fun},
there are positive constants $a$ and $b$ such that
$h_{nv}^{\overline{B}_i} \leq a h_{nv}^{\overline{B}} + b$ for all $i$,
where $\overline{B}_i$ is a polarization $(B; \overline{A}_i, \ldots, \overline{A}_i)$
given by $\overline{A}_i$.
We set
\[
 S =
\{ P \in \PP^n(\QQ(z_1, \cdots, z_d)) \mid h_{nv}^{\overline{B}}(P) \leq M \}
\]
Then, for any $P \in S$,
$h_{nv}^{\overline{B}_i}(P) \leq a M + b$.
Moreover, there are $f_0, \cdots, f_n \in \ZZ[z_1, \cdots, z_d]$
such that $f_0, \cdots, f_n$ are relatively prime and
$P = (f_0 : \cdots : f_n)$.
Here,
\[
 h_{nv}^{\overline{B}_i}(P) =
\max \{ \deg_i(f_0), \ldots, \deg_i(f_n) \} e,
\]
because $c_1(\overline{A}_i)^{\wedge d} = 0$ and
$f_0, \cdots, f_n$ are relatively prime.
Thus, there is a constant $M_1$ independent on $P$
such that $\deg_i(f_j) \leq M_1$ for
all $i,j$.
On the other hand,
\begin{multline*}
 h_{nv}^{\overline{B}}(P) =
\sum_i \max \{ \deg_i(f_0), \ldots, \deg_i(f_n) \} 
\adeg\left( \acherncl_1(\rest{\overline{H}}{\Delta^{(i)}_{\infty}})^d
\right) \\
+ \int_{(\PP^1)^d} \log \left( \max_i \{ |f_i| \} \right) c_1(\overline{H})^{\wedge d}.
\end{multline*}
Hence, there is a constant $M_2$ independent on $P$ such that
\[
 \int_{(\PP^1)^d} \log (|f_i|) c_1(\overline{H})^{\wedge d} \leq M_2
\]
for all $i$.
Thus, by Lemma~\ref{lem:estimate:value:integral},
we have only finitely many $f_i$'s as above.
\QED

\begin{Theorem}
\label{thm:northcott:thm:fun:field}
Let $K$ be a finitely generated field over $\QQ$ with $\trdeg_{\QQ}(K) = d$,
and $ \overline{B} = (B; \overline{H}_1, \ldots, \overline{H}_d)$ 
a big polarization of $K$, i.e.,
$\overline{H}_i$'s are nef and big.
Let $X$ be a projective variety over $K$, and $L$ an ample line bundle
on $X$.
Then, for any numbers $M$ and any positive integers $e$, the set
\[
\{ P \in X(\overline{K}) \mid h^{\overline{B}}_L(P) \leq M,
\quad [K(P) : K] \leq e \} 
\]
is finite.
\end{Theorem}

\Proof
We set $B_0 = (\PP_{\ZZ}^1)^d$ and
$\overline{H}_0 = p_1^*((\OO(1), \normabb_{FS})) \otimes \cdots \otimes
p_d^*((\OO(1), \normabb_{FS})))$.

\begin{Claim}
\label{claim:thm:northcott:thm:fun:field:1}
If $B = B_0$, $\overline{H}_0 = \overline{H}_1 = \cdots = \overline{H}_d$
and $e = 1$,
then the theorem holds.
\end{Claim}

If $m$ is sufficiently large, then we have an embedding
$X \hookrightarrow \PP^n$ by $L^{\otimes m}$. Thus, we may assume
$(X, L) = (\PP^n_{\QQ(z_1, \ldots, z_d)}, \OO(1))$. 
Hence, this claim follows from
Lemma~\ref{lem:northcott:special:case}.

\begin{Claim}
\label{claim:thm:northcott:thm:fun:field:2}
To prove the theorem, we may assume that $B = B_0$ and $\overline{H}_0 =
\overline{H}_1 = \cdots = \overline{H}_d$.
\end{Claim}

As in the previous claim, we may assume $(X, L) = (\PP^n_K, \OO(1))$.
Since $\trdeg_{\QQ}(K) = d$,
$K$ contains $\QQ(z_1, \ldots, z_d)$, which means that
there is a rational map $B \dasharrow B_0$. 
Thus, replacing $B$ by the graph of
$B \dasharrow B_0$, 
we may assume that there is a generically
finite morphism $\mu : B \to B_0$.
Then, since $\overline{H}_i$'s are big,
by Proposition~\ref{prop:equiv:cond:big},
there are positive integers $n_1, \ldots, n_d$ such that
$\overline{H}_i^{\otimes n_i} \otimes \mu^*(\overline{H}_0)^{\otimes -1}$
is effective for every $i$.
Thus, using (5) of Proposition~\ref{prop:elem:prop:height:fun},
$h_L^{(B; \mu^*(\overline{H}_0), \ldots, \mu^*(\overline{H}_0))}
\leq a h_L^{(B; \overline{H}_1, \ldots, \overline{H}_d)} + b$
for some positive constants $a$ and $b$.
Hence, we may assume that $\overline{H}_1 = \cdots = \overline{H}_d 
= \mu^*(\overline{H}_0)$.
Therefore, our claim follows from
Proposition~\ref{prop:base:change:field}.

\begin{Claim}
\label{claim:thm:northcott:thm:fun:field:3}
In order to prove our theorem, we may assume that $e=1$.
\end{Claim}

It is sufficient to show that the set 
\[
\{ P \in X(\overline{K}) \mid
h^{\overline{B}}_{L}(P) \leq M, \quad
[K(P) : K] = e \}
\]
is finite for any numbers $M$ and any integers $e \geq 1$.
Let $(\XX, \LL)$ be a $C^{\infty}$-model of $(X, L)$.
Let $P$ be a point of $X(\overline{K})$ with
$h^{\overline{B}}_{L}(P) \leq M$ and
$[K(P) : K] = e$.
Let $\{ \sigma_1, \ldots, \sigma_{e} \}$ be the set of all embeddings
of $K(P)$ into $\overline{K}$.
Let $P_i \in X(\overline{K})$ be a point given by
$\Spec(\overline{K}) \overset{\sigma_i^*}{\longrightarrow} \Spec(K(P))
\overset{P}{\longrightarrow} X$, and let
$\Delta_{P_i}$ be the closure of $P_i$ in $\XX$.
Then, we have $\Delta_P = \Delta_{P_i}$ for all $i$.
Let $\mathcal{Y}$ be the main component of 
${\displaystyle 
\underbrace{\XX \times_B \cdots \times_B \XX}_{\text{$e$-times}}}$, and 
$\overline{\mathcal{M}} = p_1^*(\overline{\LL})
\otimes \cdots \otimes p_e^*(\overline{\LL})$,
where $p_i$ is the projection to the $i$-th factor.
Moreover, let $f : \XX \to B$ and $f' : \mathcal{Y} \to B$ be
the canonical morphism. Then, using the projection formula,
\begin{align*}
h^{\overline{B}}_{(\mathcal{Y}, \overline{\mathcal{M}})}(P_1, \ldots,P_e)
& = \frac{
\adeg \left(
(f')^*(\acherncl_1(\overline{H}_1) \cdots \acherncl_1(\overline{H}_d)) \cdot 
\acherncl_1(\overline{\mathcal{M}}) \cdot 
(\Delta_{(P_1, \ldots, P_e)}, 0)
\right)}
{\deg(\Delta_{(P_1, \ldots, P_e)} \to B)} \\
& = \frac{
\sum_{i=1}^e
\adeg \left(
p_i^* f^*(\acherncl_1(\overline{H}_1) \cdots \acherncl_1(\overline{H}_d)) \cdot 
p_i^*\acherncl_1(\overline{\LL}) \cdot 
(\Delta_{(P_1, \ldots, P_e)}, 0) 
\right)}
{\deg(\Delta_{(P_1, \ldots, P_e)} \to B)} \\
& =
\frac{
\sum_{i=1}^e
\deg(\Delta_{(P_1, \ldots, P_e)} \to \Delta_{P_i})\adeg \left(
f^*(\acherncl_1(\overline{H}_1) \cdots \acherncl_1(\overline{H}_d)) \cdot 
\acherncl_1(\overline{\LL}) \cdot 
(\Delta_{P_i}, 0) 
\right)}
{\deg(\Delta_{(P_1, \ldots, P_e)} \to B)} \\
& = e h^{\overline{B}}_{(\XX, \LL)}(P).
\end{align*}
Let $\rho : X^e \to Z = X^e/\mathfrak{S}_e$ 
be the quotient of $X^e$ by the symmetric group $\mathfrak{S}_e$.
Since $M = \mathcal{M}_{K}$ is invariant under the action of
$\mathfrak{S}_e$, there is a line bundle $N$ on $Z$ with
$\rho^*(N) = M$.
Thus, $\rho^*(h_N^{\overline{B}}) = 
h^{\overline{B}}_{(\mathcal{Y}, \overline{\mathcal{M}})} + O(1)$.
Hence,
\[
h_N^{\overline{B}}(\rho(P_1, \ldots, P_e)) \leq
e h_L^{\overline{B}}(P) + C \leq e M + C
\]
for some constant $C$ independent on $P$.
Moreover, $\rho(P_1, \ldots, P_e)$ is defined over $K$.
Thus, if our theorem holds for the case $e=1$,
there are finitely many 
$\rho(P_1, \ldots, P_e)$ with 
\[
 h_N^{\overline{B}}(\rho(P_1, \ldots, P_e)) \leq e M + C.
\]
Here the number of the fiber of $\rho$ is $e!$ at most.
Hence, we have our claim.

\medskip
Let us start the proof of the theorem.
First, by Claim~\ref{claim:thm:northcott:thm:fun:field:2},
we may assume $B = B_0$ and $\overline{H}_0 = \overline{H}_1 = \cdots = \overline{H}_d$.
Thus, Claim~\ref{claim:thm:northcott:thm:fun:field:1} and
Claim~\ref{claim:thm:northcott:thm:fun:field:3} implies our
theorem. 
\QED

\section{Estimate of height functions in terms of intersection numbers}
Let $K$ be a finitely generated field over $\QQ$ with $d = \trdeg_{\QQ}(K)$,
$B$ an arithmetic model of $K$, and let $\overline{H}$
be a nef $C^{\infty}$ hermitian $\QQ$-line bundles on $B$.
Let $\overline{B} = (B; \overline{H}, \ldots, \overline{H})$ be a polarization
of $K$ given by $\overline{H}$.
Let $X$ be an $e$-dimensional projective variety over $K$, and
$L$ a line bundle on $X$.
Let $(\XX, \overline{\LL})$ be a $C^{\infty}$-model of $(X, L)$, and
$\pi : \XX \to B$ the canonical morphism.
The purpose of this section is to prove the following theorem.

\begin{Theorem}
\label{thm:estimate:liminf:inf:height}
We assume that $\adeg(\acherncl_1(\overline{H})^{d+1}) = 0$,
and that, for some rational number $a$,
$\overline{\LL} + a \pi^*(\overline{H})$ is vertically nef 
and $(\LL + a \pi^*(H))_{\QQ}$ is ample on $\XX_{\QQ}$.
Then, we have the following.
\begin{enumerate}
\renewcommand{\labelenumi}{(\arabic{enumi})}
\item
If $\adeg(\acherncl_1(\overline{\LL})^{e+1} \cdot 
\acherncl_1(\pi^*(\overline{H}))^d) > 0$,
then
\[
\sup_{Y \subsetneq X} \left\{ \inf_{x \in (X \setminus Y)(\overline{K})}
h_{(\XX, \overline{\LL})}^{\overline{B}}(x) \right\} \geq 0,
\]
where $Y$ runs over all proper closed subsets of $X$.

\item
If $\deg(H_{\QQ}^d) > 0$ and ${\displaystyle \inf_{x \in X(\overline{K})} 
h_{(\XX, \overline{\LL})}^{\overline{B}}(x) \geq 0}$,
then
\[
\adeg(\acherncl_1(\overline{\LL})^{e+1} \cdot 
\acherncl_1(\pi^*(\overline{H}))^d) \geq 0.
\]
\end{enumerate}
\end{Theorem}

\Proof
(1)
Since $\adeg(\acherncl_1(\overline{H})^{d+1}) = 0$, 
by virtue of the projection formula
(cf. (1) and (2) of 
Proposition~\ref{prop:projection:formula:intersection:num}),
we can easily see that
\[
h_{(\XX, \overline{\LL})}^{\overline{B}} = 
h_{(\XX, \overline{\LL} +m \pi^*(\overline{H}))}^{\overline{B}}
\]
and
\[
\adeg(\acherncl_1(\overline{\LL})^{e+1} \cdot \acherncl_1(\pi^*(\overline{H}))^d) =
\adeg\left( 
(\acherncl_1(\overline{\LL} + m \pi^*(\overline{H}))^{e+1} \cdot 
\acherncl_1(\pi^*(\overline{H}))^d
\right).
\]
Thus, we may assume that $\overline{\LL}$ is vertically nef and
$\LL_{\QQ}$ is ample on $B_{\QQ}$
because
\[
\overline{\LL} + m \pi^*(\overline{H}) = 
\overline{\LL} + a \pi^*(\overline{H}) +
(m-a) \pi^*(\overline{H}).
\]
Moreover, since
\[
\adeg \left(
\left( \acherncl_1(\overline{\LL}) + 
m \acherncl_1(\pi^*(\overline{H})) \right)^{e+d+1}
\right) = 
\binom{e+d+1}{d} \adeg \left( \acherncl_1(\overline{\LL})^{e+1} \cdot 
\acherncl_1(\pi^*(\overline{H}))^d \right) m^d + O(m^{d-1}),
\]
we may further assume that $\adeg(\acherncl_1(\overline{\LL})^{e+d+1}) > 0$.
Thus, by Theorem~\ref{thm:existense:small:sec}, 
for a sufficiently large integer $n$, there is a non-zero section
$s$ of $H^0(\XX, \LL^{\otimes n})$ with $\Vert s \Vert_{\sup} < 1$.
We set $Y = (\zero(s))_K$. Then, for any $x \in (X \setminus Y)(\overline{K})$,
\[
h_{(\XX, \overline{\LL})}^{\overline{B}}(x) =
\frac{1}{n \deg(\Delta_x \to B)} \adeg \left(
\rest{\widehat{\zero(s)}}{\Delta_x} \cdot 
\acherncl_1(\rest{\pi^*(\overline{H})}{\Delta_x})^d
\right).
\]
Here, $\rest{\widehat{\zero(s)}}{\Delta_x}$ is effective, and
$\rest{\pi^*(\overline{H})}{\Delta_x}$ is nef. Thus, 
by (2) of Proposition~\ref{prop:intersection:nef:line:bundle}, we have
$h_{(\XX, \overline{\LL})}^{\overline{B}}(x) \geq 0$ for all 
$x \in (X \setminus Y)(\overline{K})$.

\bigskip
(2)
The proof of (2) is very similar to \cite[Lemma~5.4]{ZhPos}, or
Lemma~\ref{lem:intersection:ample:nef}.
First of all, we need the following two claims.

\begin{Claim}
\label{claim:thm:estimate:liminf:inf:height:1}
Let $\nu : \mathcal{Y} \to B$ be a surjective morphism
of projective arithmetic varieties, and
let $\overline{\LL}_1$ and $\overline{\LL}_2$ be $C^{\infty}$-hermitian
$\QQ$-line bundles on $\mathcal{Y}$ such that
$\overline{\LL}_1$ and $\overline{\LL}_2$ are vertically nef, and that
$(\LL_1)_{\QQ}$ and $(\LL_2)_{\QQ}$ are ample on $\mathcal{Y}_{\QQ}$.
Let us fix an integer $s$ with $0 \leq s \leq \dim(\mathcal{Y}/B)$,
where $\dim(\mathcal{Y}/B)$ is the dimension of the generic fiber
of $\mathcal{Y} \to B$.
We assume the following:
\begin{enumerate}
\renewcommand{\labelenumi}{(\alph{enumi})}
\item
$\adeg \left( \acherncl_1(\rest{\overline{\LL}_1}{\Gamma})^{s+1} \cdot
\acherncl_1(\rest{\nu^*(\overline{H})}{\Gamma})^{d} \right)
\geq 0$ for any integral subschemes $\Gamma$ on $\mathcal{Y}$
with $\nu(\Gamma) = B$ and $\dim(\Gamma/B) = s$.

\item
$\adeg \left( \acherncl_1(\rest{\overline{\LL}_2}{\Gamma})^{\dim(\Gamma/B)+1} 
\cdot
\acherncl_1(\rest{\nu^*(\overline{H})}{\Gamma})^{d} \right)
> 0$ for any integral subschemes $\Gamma$ on $\mathcal{Y}$
with $\nu(\Gamma) = B$.
\end{enumerate}
Then,
$\adeg \left( \acherncl_1(\overline{\LL}_1)^{s+1} \cdot
\acherncl_1(\overline{\LL}_2)^{\dim(\mathcal{Y}/B)-s} \cdot
\acherncl_1(\nu^*(\overline{H}))^d \right) \geq 0$.
\end{Claim}

\Proof
We prove this claim by induction on $\dim(\mathcal{Y}/B)$.
If $s = \dim(\mathcal{Y}/B)$, then our assertion is trivial.
Hence we may assume that $\dim(\mathcal{Y}/B) > s$.
Since
\[
\adeg \left( \acherncl_1(\overline{\LL}_2)^{\dim(\mathcal{Y}/B)+1} 
\cdot \acherncl_1(\nu^*(\overline{H}))^d \right) > 0,
\]
in the same way as the proof of (1),
\[
\adeg \left(
\left( \acherncl_1(\overline{\LL}_2) + 
m \acherncl_1(\nu^*(\overline{H})) \right)^{\dim(\mathcal{Y}/B)+d+1}
\right) > 0
\]
for a sufficiently large $m$.
Thus, by Theorem~\ref{thm:existense:small:sec}, 
for a sufficiently large integer $n$, there is a non-zero section
$s$ of $H^0(\XX, n(\LL_2 + m \nu^*(H))))$ 
with $\Vert s \Vert_{\sup} < 1$.
Let $\zero(s) = \sum_{i} a_i \Gamma_i$ be the irreducible decomposition
as cycles.
Since $\adeg \left( \acherncl_1(\overline{H})^{d+1} \right) = 0$,
we have
\begin{multline*}
\adeg \left( \acherncl_1(\overline{\LL}_1)^{s+1} \cdot
\acherncl_1(\overline{\LL}_2)^{\dim(\mathcal{Y}/B)-s} \cdot
\acherncl_1(\nu^*(\overline{H}))^d \right) = \\
\adeg \left( \acherncl_1(\overline{\LL}_1)^{s+1} \cdot
\acherncl_1(\overline{\LL}_2)^{\dim(\mathcal{Y}/B)-s-1} \cdot
\acherncl_1(\overline{\LL}_2 + m \nu^*(\overline{H})) \cdot
\acherncl_1(\nu^*(\overline{H}))^d \right),
\end{multline*}
which implies that
\begin{multline*}
n \adeg \left( \acherncl_1(\overline{\LL}_1)^{s+1} \cdot
\acherncl_1(\overline{\LL}_2)^{\dim(\mathcal{Y}/B)-s} \cdot
\acherncl_1(\nu^*(\overline{H}))^d \right) = \\
\sum_{i} a_i \adeg \left( 
\acherncl_1(\rest{\overline{\LL}_1}{\Gamma_i})^{s+1} \cdot
\acherncl_1(\rest{\overline{\LL}_2}{\Gamma_i})^{\dim(\mathcal{Y}/B)-s-1} \cdot
\acherncl_1(\rest{\nu^*(\overline{H})}{\Gamma_i})^d \right) \\
+ \int_{\XX(\CC)} (-\log \Vert s \Vert_{\sup})
c_1(\overline{\LL}_1)^{\wedge s+1} \wedge 
c_1(\overline{\LL}_2)^{\wedge \dim(\mathcal{Y}/B)-s-1} \wedge
c_1(\nu^*(\overline{H}))^{\wedge d} 
\end{multline*}
Thus, by Lemma~\ref{lem:non:negative:prod:semi:positive},
it is sufficient to show
\addtocounter{Claim}{1}
\begin{equation}
\label{lem:thm:estimate:liminf:inf:height:1:eqn}
\adeg \left( 
\acherncl_1(\rest{\overline{\LL}_1}{\Gamma_i})^{s+1} \cdot
\acherncl_1(\rest{\overline{\LL}_2}{\Gamma_i})^{\dim(\mathcal{Y}/B)-s-1} \cdot
\acherncl_1(\rest{\nu^*(\overline{H})}{\Gamma_i})^d \right)
\geq 0
\end{equation}
for all $i$.

If $\Gamma_i$ maps surjectively to $B$, 
then, by hypothesis of induction, we can see 
\eqref{lem:thm:estimate:liminf:inf:height:1:eqn}.
Thus, we assume that $\Gamma_i$ dose not map surjectively to $B$.
Let $T$ be the generic fiber of $\Gamma_i \to \nu(\Gamma_i)$.
Since $\dim T \geq e$, we can see
\begin{multline*}
\adeg \left( 
\acherncl_1(\rest{\overline{\LL}_1}{\Gamma_i})^{s+1} \cdot
\acherncl_1(\rest{\overline{\LL}_2}{\Gamma_i})^{\dim(\mathcal{Y}/B)-s-1} \cdot
\acherncl_1(\rest{\nu^*(\overline{H})}{\Gamma_i})^d \right) \\
= \begin{cases}
\deg(\rest{\LL_1}{T}^{s+1} \cdot \rest{\LL_2}{T}^{\dim(\mathcal{Y}/B)-s-1}) 
\adeg(\acherncl_1(\rest{\overline{H}}{\nu(\Gamma_i)})^d) 
& \text{if $\dim T = \dim(\mathcal{Y}/B)$} \\
0 & \text{if $\dim T > \dim(\mathcal{Y}/B)$}.
\end{cases}
\end{multline*}
Here $\overline{\LL}_1$ and $\overline{\LL}_2$ are vertically nef and
$\overline{H}$ is nef. Thus, by the above formula,
we have \eqref{lem:thm:estimate:liminf:inf:height:1:eqn}
even if $\Gamma_i$ dose not map surjectively to $B$.
\QED

\begin{Claim}
\label{claim:thm:estimate:liminf:inf:height:2}
If $\deg(H_{\QQ}^d) > 0$, then
there is an ample $C^{\infty}$-hermitian line bundle $\overline{\mathcal{M}}$ 
on $\XX$ such that
\[
\adeg \left( \
\acherncl_1(\rest{\overline{\mathcal{M}}}{\Gamma})^{\dim(\Gamma/B)+1} 
\cdot
\acherncl_1(\rest{\pi^*(\overline{H})}{\Gamma})^{d} \right)
> 0
\]
for any integral subschemes $\Gamma$ on $\XX$
with $\pi(\Gamma) = B$.
\end{Claim}

\Proof
Let $\overline{\mathcal{N}}$ be an ample
$C^{\infty}$-hermitian line bundle on $\XX$,
and let $\normabb$ be the metric of $\overline{H}$.
For a positive number $c$ with $0 < c < 1$, we set
$\overline{A} = (H, c \normabb)$.
Then, since $\deg(H_{\QQ}^d) > 0$, we can see that
$\adeg(\acherncl_1(\overline{A}) \cdot \acherncl_1(\overline{H})^d) > 0$.
Let $\Gamma$ be a subscheme of $\XX$ with $\pi(\Gamma) = B$.
Then, by (1) of Proposition~\ref{prop:intersection:nef:line:bundle},
\[
\adeg \left( \
\acherncl_1(\rest{\overline{\mathcal{N}}}{\Gamma})^{i} \cdot
\acherncl_1(\rest{\pi^*(\overline{A})}{\Gamma})^{\dim(\Gamma/B)+1-i} 
\cdot
\acherncl_1(\rest{\pi^*(\overline{H})}{\Gamma})^{d} \right)
\geq 0
\]
for all $0 \leq i \leq \dim(\Gamma/B)+1$.
Further,
\begin{multline*}
\adeg \left( \
\acherncl_1(\rest{\overline{\mathcal{N}}}{\Gamma})^{\dim(\Gamma/B)} \cdot
\acherncl_1(\pi^*(\rest{\overline{A})}{\Gamma}) 
\cdot
\acherncl_1(\rest{\pi^*(\overline{H})}{\Gamma})^{d} \right)
\\
= \deg((\mathcal{N}_{\eta})^{\dim(\Gamma/B)})
\adeg(\acherncl_1(\overline{A}) \cdot \acherncl_1(\overline{H})^d)
>  0,
\end{multline*}
where $\eta$ means the restriction to the generic fiber of
$\Gamma \to B$.
Thus, if we set
$\overline{\mathcal{M}} = \overline{\mathcal{N}} + \pi^*(\overline{A})$,
then we have the desired hermitian line bundle.
\QED

\bigskip
Let us go back to the proof of (2).
We prove (2) by induction on $e$.
If $e = 0$, then the assertion is trivial.
Thus, we assume $e > 0$.

In the same way as in the proof of (1), we may assume that
$\overline{\LL}$ is vertically nef and $\LL_{\QQ}$ is ample on $\XX_{\QQ}$.
By Claim~\ref{claim:thm:estimate:liminf:inf:height:2},
there is an ample $C^{\infty}$-hermitian line bundle $\overline{\mathcal{M}}$ 
on $\XX$ such that
\[
\adeg \left( \
\acherncl_1(\rest{\overline{\mathcal{M}}}{\Gamma})^{\dim(\Gamma/B)+1} 
\cdot
\acherncl_1(\rest{\pi^*(\overline{H})}{\Gamma})^{d} \right)
> 0
\]
for any integral subschemes $\Gamma$ on $\XX$
with $\pi(\Gamma) = B$.
Thus, by hypothesis of induction and
Claim~\ref{claim:thm:estimate:liminf:inf:height:1},
for any integral subschemes $\Gamma$ on $X$ with
$\pi(\Gamma) = B$ and $\dim(\Gamma/B) < e$,
\addtocounter{Claim}{1}
\begin{equation}
\label{eqn:thm:estimate:liminf:inf:height:1}
\adeg \left( \left(\acherncl_1(\rest{\overline{\LL}}{\Gamma}) + 
t \acherncl_1(
\rest{\overline{\mathcal{M}}}{\Gamma})\right)^{\dim(\Gamma/B)+1} \cdot
\acherncl_1(\rest{\pi^*(\overline{H})}{\Gamma})^d \right) > 0
\end{equation}
for all $t > 0$.
Moreover, we can see
\addtocounter{Claim}{1}
\begin{equation}
\label{eqn:thm:estimate:liminf:inf:height:2}
\adeg \left(\acherncl_1(\overline{\LL})^{e+1-i} \cdot
\acherncl_1(\overline{\mathcal{M}})^{i} \cdot
\acherncl_1(\pi^*(\overline{H}))^d \right) \geq 0
\end{equation}
for all $1 \leq i \leq e+1$.
We set
\[
P(t) = \adeg \left( \left(\acherncl_1(\overline{\LL}) + 
t \acherncl_1(\overline{\mathcal{M}})\right)^{e+1} \cdot
\acherncl_1(\pi^*(\overline{H}))^d \right)
\]
Here we claim the following.

\begin{Claim}
\label{claim:thm:estimate:liminf:inf:height:3}
If $t > 0$ and $P(t) > 0$, then
$P(t) \geq t^{e+1} \adeg \left(
\acherncl_1(\overline{\mathcal{M}})^{e+1} \cdot
\acherncl_1(\pi^*(\overline{H}))^d \right)$.
\end{Claim}

Clearly we may assume that $t$ is a rational number.
For simplicity, we set $\overline{\mathcal{N}} =
\overline{\LL} + t \overline{\mathcal{M}}$. Then,
by \eqref{eqn:thm:estimate:liminf:inf:height:1} and
$P(t) > 0$,
we have
\[
\adeg \left( \
\acherncl_1(\rest{\overline{\mathcal{N}}}{\Gamma})^{\dim(\Gamma/B)+1} 
\cdot
\acherncl_1(\rest{\pi^*(\overline{H})}{\Gamma})^{d} \right)
> 0
\]
for any integral subschemes $\Gamma$ on $\XX$
with $\pi(\Gamma) = B$.
Thus, by using Claim~\ref{claim:thm:estimate:liminf:inf:height:1} and
the assumption
${\displaystyle \inf_{x \in X(\overline{K})} 
h_{(\XX, \overline{\LL})}^{\overline{B}}(x) \geq 0}$,
we obtain $\adeg \left( \acherncl_1(\overline{\LL}) \cdot
\acherncl_1(\overline{\mathcal{N}})^e \cdot \acherncl_1(\pi^*(\overline{H}))
\right) \geq 0$.
Therefore, by using \eqref{eqn:thm:estimate:liminf:inf:height:2},
\begin{align*}
P(t) & = \adeg \left( (\acherncl_1(\overline{\LL}) + t
\acherncl_1(\overline{\mathcal{M}})) \cdot \acherncl_1(\overline{\mathcal{N}})^e \cdot
\acherncl_1(\pi^*(\overline{H}))^d \right) \\
& \geq
t \adeg \left( \acherncl_1(\overline{\mathcal{M}}) \cdot 
\acherncl_1(\overline{\mathcal{N}})^e \cdot
\acherncl_1(\pi^*(\overline{H}))^d \right) \\
& \geq t^{e+1} \adeg \left( \acherncl_1(\overline{\mathcal{M}})^{e+1} \cdot
\acherncl_1(\pi^*(\overline{H}))^d \right).
\end{align*}
Thus, we get the claim.

\medskip
Let $t_0 = \max \{ t \in \RR \mid P(t) = 0 \}$.
We assume $t_0 > 0$. Then, by the above claim, for any $t > t_0$,
\[
P(t) \geq t^{e+1} \adeg \left(
\acherncl_1(\overline{\mathcal{M}})^{e+1} \cdot
\acherncl_1(\pi^*(\overline{H}))^d \right).
\]
Thus, taking $t \to t_0$,
\[
0 = P(t_0) \geq t_0^{e+1} \adeg \left(
\acherncl_1(\overline{\mathcal{M}})^{e+1} \cdot
\acherncl_1(\pi^*(\overline{H}))^d \right) > 0.
\]
This is a contradiction. Therefore, $t_0 \leq 0$.
In particular, $P(0) \geq 0$,
which is nothing more than the assertion of (2).
\QED

As corollary, we have the following generalization of
\cite[Theorem~(5.2)]{ZhPos}.

\begin{Corollary}
\label{cor:estimate:liminf:inf:height}
We assume that
$\adeg(\acherncl_1(\overline{H})^{d+1}) = 0$,
$\deg(H_{\QQ}^d) > 0$, and that,
for some rational number $a$,
$\overline{\LL} + a \pi^*(\overline{H})$ is vertically nef 
and $(\LL + a \pi^*(H))_{\QQ}$ is ample on $\XX_{\QQ}$.
Then we have the following inequalities:
\[
\sup_{Y \subsetneq X} \left\{ \inf_{x \in (X \setminus Y)(\overline{K})}
h_{(\XX, \overline{\LL})}^{\overline{B}}(x) \right\} \geq
\frac{\adeg(\acherncl_1(\overline{\LL})^{e+1} \cdot 
\acherncl_1(\pi^*(\overline{H}))^d)}
{(e+1)\deg(\LL_K^e)}
\geq
\inf_{x \in X(\overline{K})} h_{(\XX, \overline{\LL})}^{\overline{B}}(x).
\]
\end{Corollary}

\Proof
Let $c$ be a real number with $0 < c < 1$.
We set $\overline{A} = (H, c\normabb)$, where
$\normabb$ is the metric of $\overline{H}$.
Then, 
\[
 \adeg (\acherncl_1(\overline{A}) \cdot \acherncl_1(\overline{H})^d) > 0
\]
because $\deg(H_{\QQ}^d) > 0$.
Let $\lambda$ be an arbitrary rational number with
\[
\lambda < 
\frac{\adeg(\acherncl_1(\overline{\LL})^{e+1} \cdot 
\acherncl_1(\pi^*(\overline{H}))^d)}
{(e+1)\deg(\LL_K^e) \adeg (\acherncl_1(\overline{A}) \cdot 
\acherncl_1(\overline{H})^d)}.
\]
Then, it is easy to see that
\[
\adeg\left(
(\acherncl_1(\overline{\LL})- \lambda \acherncl_1(\pi^*(\overline{A})) )^{e+1}
\cdot \acherncl_1(\pi^*(\overline{H}))^d \right) > 0.
\]
Here, note that
\[
\overline{\LL} - \lambda \pi^*(\overline{A}) + 
(a + \lambda) \pi^*(\overline{H}) =
\overline{\LL} +
\lambda \pi^*(\overline{H} - \overline{A}) + a \pi^*(\overline{H})
\]
is vertically nef and  ample on $\XX_{\QQ}$
because $c_1(\overline{H}) = c_1(\overline{A})$.
Thus, applying (1) of Theorem~\ref{thm:estimate:liminf:inf:height},
\[
\sup_{Y \subsetneq X} \left\{ \inf_{x \in (X \setminus Y)(\overline{K})}
h_{(\XX, \overline{\LL} - \lambda \pi^*(\overline{A}))}^{\overline{B}}(x)
\right\} \geq 0,
\]
which implies
\[
\sup_{Y \subsetneq X} \left\{ \inf_{x \in (X \setminus Y)(\overline{K})}
h_{(\XX, \overline{\LL})}^{\overline{B}}(x) \right\} \geq
\lambda \adeg (\acherncl_1(\overline{A}) \cdot \acherncl_1(\overline{H})^d).
\]
Thus, we get
\[
\sup_{Y \subsetneq X} \left\{ \inf_{x \in (X \setminus Y)(\overline{K})}
h_{(\XX, \overline{\LL})}^{\overline{B}}(x) \right\} \geq
\frac{\adeg(\acherncl_1(\overline{\LL})^{e+1} \cdot 
\acherncl_1(\pi^*(\overline{H}))^d)}
{(e+1)\deg(\LL_K^e)}.
\]

\medskip
Next let $\mu$ be an arbitrary rational number with 
\[
\mu \leq 
\frac{\inf_{x \in X(\overline{K})} h_{(\XX, \overline{\LL})}^{\overline{B}}(x)}
{\adeg (\acherncl_1(\overline{A}) \cdot \acherncl_1(\overline{H})^d)}.
\]
Then,
\[
\inf_{x \in X(\overline{K})}
h_{(\XX, \overline{\LL} - \mu \pi^*(\overline{A}) )}^{\overline{B}}(x) \geq 0.
\]
Thus, by (2) of Theorem~\ref{thm:estimate:liminf:inf:height},
\[
\adeg\left(
(\acherncl_1(\overline{\LL})- \mu \acherncl_1(\pi^*(\overline{A})) )^{e+1}
\cdot \acherncl_1(\pi^*(\overline{H}))^d \right) \geq 0,
\]
which says us that
\[
\adeg\left(
\acherncl_1(\overline{\LL})^{e+1}
\cdot \acherncl_1(\pi^*(\overline{H}))^d \right) \geq
\mu(e+1)\deg(\LL_K^{e}) \adeg (\acherncl_1(\overline{A}) \cdot 
\acherncl_1(\overline{H})^d).
\]
Hence, we get the second inequality.
\QED

\section{Equidistribution theorem over finitely generated fields}
Let $K$ be a finitely generated field over $\QQ$ with $d = \trdeg_{\QQ}(K)$,
$B$ an arithmetic model of $K$, and let $\overline{H}$
be a nef $C^{\infty}$ hermitian $\QQ$-line bundles on $B$.
Let $\overline{B} = (B; \overline{H}, \ldots, \overline{H})$ be a polarization
of $K$ given by $\overline{H}$.
Let $X$ be an $e$-dimensional projective variety over $K$.
Let $\{ x_m \}_{m=1}^{\infty}$ be a sequence of
elements of $X(\overline{K})$.
We say $\{ x_m \}$ is {\em generic} if
any subsequences of $\{ x_m \}$ are not contained
in any proper closed subsets of $X(\overline{K})$.

Let $L$ be a line bundle on $X$.
Let $(\XX, \overline{\LL})$ be a $C^{\infty}$-model of $(X, L)$, and
$\pi : \XX \to B$ the canonical morphism.
Then, we have the following equidistribution theorem,
which is a generalization of Szpiro-Ullmo-Zhang's result
(cf. \cite{SUZEqui}, \cite{UlPos} and \cite{ZhEqui}).

\begin{Theorem}
\label{thm:equi:dist}
Let $h : X(\overline{K}) \to \RR$ be a representative of
the class of height functions associated with
$L$ and $\overline{B}$, and 
let $(\XX_n, \overline{\LL}_n)$ be a sequence of $C^{\infty}$-models
of $(X, L)$ over $B$.
We assume the following.
\begin{enumerate}
\renewcommand{\labelenumi}{(\arabic{enumi})}
\item
$\adeg(\acherncl_1(\overline{H})^{d+1}) = 0$
and $\deg(H_{\QQ}^d) > 0$.

\item
$h(x) \geq 0$ for all $x \in X(\overline{K})$.

\item
There is a Zariski open set $U$ of $B$ such that
$(\XX_n)_U = \XX_U$ for all $n$.

\item
$\sup_{x \in X(\overline{K})} 
\vert h(x) - h_{(\XX_n, \overline{\LL}_n)}^{\overline{B}}(x) \vert$
converges to $0$ as $n$ tends to $\infty$.

\item
For $n \gg 0$,
$\overline{\LL}_n$
is vertically nef, and $(\LL_n)_{\QQ}$ is ample on $(\XX_n)_{\QQ}$.

\item
There are a connected open set $W$ of $U(\CC)$
\rom{(}in the topology as analytic spaces\rom{)}
and a positive $C^{\infty}$-form $\omega$ on $\pi^{-1}(W)$ such that
$c_1(\overline{\LL}_n) = \omega$ on $\pi^{-1}(W)$ for $n \gg 0$.
\end{enumerate}
Let $\{ x_m \}$ be a generic sequence in $X(\overline{K})$ with
$\lim_{m \to \infty} h(x_m) = 0$.
Then, over $\pi^{-1}(W)$, we have the following weak convergence
\[
\lim_{m \to \infty}
\frac{\delta_{\Delta_{x_m}} \wedge \pi^*(c_1(\overline{H}))^{\wedge d}}
{\deg(\Delta_{x_m} \to B)}
= \left[ \frac{\omega^{\wedge e} \wedge \pi^*(c_1(\overline{H}))^{\wedge d}}{\deg(L^e)}
\right]
\]
as currents. 
\end{Theorem}

\Proof
Let $f$ be a real valued $C^{\infty}$-function on $\pi^{-1}(W)$ with compact support.
We need to show that
\[
\lim_{m \to \infty} \frac{\int_{\pi^{-1}(W)} f \delta_{\Delta_{x_m}} \wedge 
\pi^*(c_1(\overline{H}))^{\wedge d}}{\deg(\Delta_{x_m} \to B)}
= \frac{\int_{\pi^{-1}(W)} f \omega^{\wedge e} \wedge 
\pi^*(c_1(\overline{H}))^{\wedge d}}{\deg(L^e)}.
\]
Let $F_{\infty}$ be the Frobenius map given by the complex conjugation.
We set $W' = W \cup F_{\infty}(W)$. Then, since
$\overline{\LL}_n$ is invariant under $F_{\infty}$,
there is a positive form $\omega'$ on $\pi^{-1}(W')$ with
$c_1(\overline{\LL}_n) = \omega'$ on $\pi^{-1}(W')$ for $n \gg 0$.
Moreover, since $\Delta_{x_m}$, $\omega'$ and $c_1(\overline{H})$
are compatible with the action induced by $F_{\infty}$, we can see that if we set
${\displaystyle f' = \frac{f + F_{\infty}^*(f)}{2}}$
(which is invariant under $F_{\infty}$), then
\[
\int_{\pi^{-1}(W)} f \delta_{\Delta_{x_m}} \wedge 
\pi^*(c_1(\overline{H}))^{\wedge d} =
\int_{\pi^{-1}(W')} f' \delta_{\Delta_{x_m}} \wedge 
\pi^*(c_1(\overline{H}))^{\wedge d}
\]
and
\[
\int_{\pi^{-1}(W)} f \omega^{\wedge e} \wedge 
\pi^*(c_1(\overline{H}))^d =
\int_{\pi^{-1}(W')} f' {\omega'}^{\wedge e} \wedge 
\pi^*(c_1(\overline{H}))^{\wedge d}.
\]
Thus, it is sufficient to see that
\[
\lim_{m \to \infty} \frac{\int_{\pi^{-1}(W')} f' \delta_{\Delta_{x_m}} \wedge 
\pi^*(c_1(\overline{H}))^{\wedge d}}{\deg(\Delta_{x_m} \to B)}
= \frac{\int_{\pi^{-1}(W')} f' {\omega'}^{\wedge e} \wedge 
\pi^*(c_1(\overline{H}))^{\wedge d}}{\deg(L^e)}.
\]

First of all, there is a positive number $\lambda_0$ such that
for all $\lambda$ with $|\lambda| \leq \lambda_0$,
$\omega' + \lambda dd^c(f')$ is semipositive on $\pi^{-1}(W')$.
Let $\OO_n(\lambda f')$ be the hermitian line bundle for $\OO_{\XX_n}$ such that
the length of $1$ at each point is given by $\exp(-\lambda f')$.
Note that since the closure of $\{ x \in \pi^{-1}(W') \mid f'(x) \not = 0 \}$
is contained in $\pi^{-1}(W') \subseteq \XX_n(\CC)$,
we may view $f'$ as a $C^{\infty}$-function on $\XX_n(\CC)$.
Here, we set 
$\overline{\LL}^{\lambda}_n = \overline{\LL}_n \otimes \OO_n(\lambda f')$.
Then, by our construction, $\overline{\LL}^{\lambda}_n$ is vertically nef and
ample on $(\XX_n)_{\QQ}$ for $n \gg 0$.
Moreover,
\[
h_{(\XX_n, \overline{\LL}^{\lambda}_n)}^{\overline{B}}(x) =
h_{(\XX_n, \overline{\LL}_n)}^{\overline{B}}(x) + 
\frac{\lambda}{\deg(\Delta_x \to B)}
\int_{\pi^{-1}(W')} f' \delta_{\Delta_x} \wedge \pi^*(c_1(\overline{H}))^{\wedge d}
\]
and
\begin{multline*}
\frac{\adeg(\acherncl_1(\overline{\LL}^{\lambda}_n)^{e+1} \cdot 
\acherncl_1(\pi^*(\overline{H}))^d)}
{(e+1)\deg(L^e)} =
\frac{\adeg(\acherncl_1(\overline{\LL}_n)^{e+1} \cdot 
\acherncl_1(\pi^*(\overline{H}))^d)}
{(e+1)\deg(L^e)} + \\
\frac{\lambda}{\deg(L^e)} \int_{\pi^{-1}(W')} f' {\omega'}^{\wedge e} \wedge 
\pi^*(c_1(\overline{H}))^{\wedge d} + O(\lambda^2),
\end{multline*}
where $\Delta_x$ is the closure of $x$ in $\XX_n$, and 
the term $O(\lambda^2)$ is
independent on $n$.
Let $\epsilon$ be an arbitrary positive number.
Then, there is a positive integer $n_1$ such that, for all $n \geq n_1$,
\[
h - \epsilon \leq h_{(\XX_n, \overline{\LL}_n)}^{\overline{B}} \leq 
h + \epsilon.
\]
On the other hand, since $\{ x_m \}$ is generic, if $n \gg 0$, 
by Corollary~\ref{cor:estimate:liminf:inf:height},
we can see
\begin{align*}
\liminf_{m} h^{\overline{B}}_{(\XX_n, \overline{\LL}^{\lambda}_n)}(x_m)
& \geq  \frac{\adeg(\acherncl_1(\overline{\LL}^{\lambda}_n)^{e+1} \cdot 
\acherncl_1(\pi^*(\overline{H}))^d)}
{(e+1)\deg(L^e)} \\
& \geq
\inf_{x \in X(\overline{K})} h^{\overline{B}}_{(\XX,\overline{\LL}_n)}(x) 
+
\frac{\lambda}{\deg(L^e)} 
\int_{\pi^{-1}(W')} f' {\omega'}^{\wedge e} \wedge \pi^*(c_1(\overline{H}))^{\wedge d} + O(\lambda^2),
\end{align*}
which implies
\[
\epsilon + 
\lambda \liminf_m \frac{\int_{\pi^{-1}(W')} f' \delta_{\Delta_{x_m}} \wedge 
\pi^*(c_1(\overline{H}))^{\wedge d}}{\deg(\Delta_{x_m} \to B)}
\geq -\epsilon + 
\lambda \frac{\int_{\pi^{-1}(W')} f' {\omega'}^{\wedge e} \wedge 
\pi^*(c_1(\overline{H}))^{\wedge d}}{\deg(L^e)}  + O(\lambda^2).
\]
Thus,
\[
\lambda \liminf_m \frac{\int_{\pi^{-1}(W')} f' \delta_{\Delta_{x_m}} \wedge 
\pi^*(c_1(\overline{H}))^{\wedge d}}{\deg(\Delta_{x_m} \to B)}
\geq
\lambda \frac{\int_{\pi^{-1}(W')} f' {\omega'}^{\wedge e} \wedge 
\pi^*(c_1(\overline{H}))^{\wedge d}}{\deg(L^e)}  + O(\lambda^2).
\]
Therefore, taking $\lambda \to 0$, we obtain
\[
\liminf_m \frac{\int_{\pi^{-1}(W')} f' \delta_{\Delta_{x_m}} \wedge 
\pi^*(c_1(\overline{H}))^{\wedge d}}{\deg(\Delta_{x_m} \to B)}
\geq \frac{\int_{\pi^{-1}(W')} f' {\omega'}^{\wedge e} \wedge 
\pi^*(c_1(\overline{H}))^{\wedge d}}{\deg(L^e)} .
\]
The above inequality still holds even if we replace $f'$ by $-f'$.
Thus,
\[
\limsup_m \frac{\int_{\pi^{-1}(W')} f' \delta_{\Delta_{x_m}} \wedge 
\pi^*(c_1(\overline{H}))^{\wedge d}}{\deg(\Delta_{x_m} \to B)}
\leq \frac{\int_{\pi^{-1}(W')} f' {\omega'}^{\wedge e} \wedge 
\pi^*(c_1(\overline{H}))^{\wedge d}}{\deg(L^e)} .
\]
Therefore, 
\[
\lim_{m\to\infty} \frac{\int_{\pi^{-1}(W')} f' \delta_{\Delta_{x_m}} \wedge 
\pi^*(c_1(\overline{H}))^{\wedge d}}{\deg(\Delta_{x_m} \to B)}
= \frac{\int_{\pi^{-1}(W')} f' {\omega'}^{\wedge e} \wedge 
\pi^*(c_1(\overline{H}))^{\wedge d}}{\deg(L^e)}.
\]
\QED

\section{Construction of the canonical height in terms
of Arakelov Geometry}
Let $K$ be a finitely generated field over $\QQ$ with $d = \trdeg_{\QQ}(K)$,
$B$ an arithmetic model of $K$, and let $\overline{H}$
be a nef $C^{\infty}$ hermitian $\QQ$-line bundles on $B$.
Let $\overline{B} = (B; \overline{H}, \ldots, \overline{H})$ be a polarization
of $K$ given by $\overline{H}$.
We assume that $\adeg(\acherncl_1(\overline{H})^{d+1}) = 0$ and
$\deg(H_{\QQ}^d) > 0$. 
Let $A$ be an abelian variety over $K$ of dimension $g$.
Let us fix a projective embedding $\iota : A \hookrightarrow \PP^N_K$,
so that we have a new embedding 
$\iota' : A \to \PP^N_K \times_K \PP^N_K$ given by
$\iota'(x) = (\iota(x), \iota([-1](x)))$.
Then, $L = {\iota'}^*(p_1^*(\OO(1)) \otimes p_2^*(\OO(1)))$ 
is ample and symmetric, where $p_i$ is the projection to
the $i$-th factor.
In this section, we would like to show the following proposition.

\begin{Proposition}
\label{prop:construction:good:seq:model}
There is a sequence of $C^{\infty}$-models $(\mathcal{A}_n, \overline{\LL}_n)$
of $(A, L)$ with the following properties:
\begin{enumerate}
\renewcommand{\labelenumi}{(\arabic{enumi})}
\item
There is a Zariski open set $U$ of $B$ such that
$(\mathcal{A}_n)_U = (\mathcal{A}_1)_U$ for all $n$, and that
$(\mathcal{A}_1)_U \to U$ is an abelian scheme over $U$.

\item
If $n$ is sufficiently large, then
$\LL_n$ is ample and
$\overline{\LL}_n$ is vertically nef.

\item
${\displaystyle \lim_{n \to \infty} 
\sup_{x \in A(\overline{K})} 
\vert \hat{h}^{\overline{B}}_L(x) - 
h^{\overline{B}}_{(\mathcal{A}_n, \overline{\LL}_n)}(x) \vert
= 0}$.

\item
There are a connected open set $W$ of $U(\CC)$,
and a positive $C^{\infty}$-form $\omega$ on $(\mathcal{A}_1)_W$ such that
$W$ is non-singular, $c_1(\overline{H})$ is positive on $W$, and
that $c_1(\overline{\LL}_n) = \omega$ on 
$(\mathcal{A}_1)_W$ for all $n \gg 0$.
\end{enumerate}
\end{Proposition}

\Proof
Let $\mathcal{A}$ be the closure of $\iota'(A)$ in
$\PP^N_B \times_B \PP^N_B$, and
$\LL = \rest{p_1^*(\OO(1)) \otimes p_2^*(\OO(1))}{\mathcal{A}}$.
Then, $\mathcal{A}_K = A$, $L = \LL_K$ is ample and
symmetric, and $\LL$ is $\pi$-ample, where
$\pi : \mathcal{A} \to B$ is the canonical morphism.
Replacing $\LL$ by $\LL \otimes \pi^*(Q)$ for some ample
line bundle $Q$ on $B$,
we may assume that $\LL$ is ample on $\mathcal{A}$.
Let $U$ be a Zariski open set of $B$ such that
$\mathcal{A}_U \to U$ is an abelian scheme and
$[-1]^*(\LL_U) = \LL_U$ over $\mathcal{A}_U$.
Shrinking $U$ if necessarily,
we may assume that $[2]^*(\LL_U) = \LL_U^{\otimes 4}$.
Let $\normabb_0$ be a hermitian metric of $\LL$ such that
$c_1(\LL, \normabb_0)$ is positive on $\mathcal{A}(\CC)$.
We would like to slightly change $\normabb_0$ on $\mathcal{A}(\CC)$.
Let $F_{\infty} : \mathcal{A}(\CC) \to \mathcal{A}(\CC)$ be the
Frobenius map given by the complex conjugation.
Since $\deg(H_{\QQ}^d) = \int_{B(\CC)} c_1(\overline{H})^{\wedge d} > 0$ and
$c_1(\overline{H})$ is semipositive, we can find
a small open set $W_1$ of $U(\CC)$ in the classical topology
such that $W_1$ is non-singular, 
$c_1(\overline{H})$ is positive on $W_1$ 
and $W_1 \cap F_{\infty}(W_1) = \emptyset$.
Here, we give a $C^{\infty}$-family of cubic metrics 
$\normabb'_{cub}$ of $\LL_{W_1}$ over $W_1$.
Then, there is a positive $C^{\infty}$-function $\rho$ on $W_1$
with $[2]^*(\LL_{W_1}, \normabb'_{cub}) = 
(\LL_{W_1}^{\otimes 4}, \rho (\normabb'_{cub})^{\otimes 4})$.
If we set $\normabb_{cub} = \rho^{1/3} \normabb'_{cub}$,
then $[2]^*(\LL_{W_1}, \normabb_{cub}) = 
(\LL_{W_1}, \normabb_{cub})^{\otimes 4}$.
Let us choose open sets $W_3$ and $W_2$ with
$W_3 \Subset W_2 \Subset W_1$.
We also choose a $C^{\infty}$ function $\phi$ on $B(\CC)$
such that $0 \leq \phi \leq 1$, $\phi \equiv 1$ on $W_3$ and
$\phi \equiv 0$ on $B(\CC) \setminus W_2$.
Let $a$ be a positive $C^{\infty}$-function on $\pi^{-1}(W_1)$
given by the equation
$\normabb_{cub} = a \normabb_0$. Here we set 
$\normabb_1 = a^{\pi^*(\phi)}\normabb_0$.
Then, $\normabb_1$ gives rise to a $C^{\infty}$-metric of $\LL$,
which coincides with $\normabb_0$ on $\pi^{-1}(B(\CC) \setminus W_2)$.
Here we claim the following.

\begin{Claim}
\label{claim:positive:limit}
For any $\epsilon > 0$, there is a positive integer $n_0$ such that,
for all $n \geq n_0$,
\[
2^{-2n}([2]^n)^*(c_1(\LL, \normabb_1)) + 
\epsilon \pi^*(c_1(\overline{H}))
\]
is positive on $\mathcal{A}_{W_2}$.
\end{Claim}

Note that the relative tangent bundle $T_{\mathcal{A}_{W_2}/W_2}$ 
is a vector subbundle of
the tangent bundle $T_{\mathcal{A}_{W_2}}$.
Let $\omega$ and $\omega'$ be the restriction of 
$c_1(\LL, \normabb_{cub})$
and $c_1(\LL, \normabb_0)$ to $T_{\mathcal{A}_{W_2}/W_2}$.
Then, $\omega$ and $\omega'$ are positive hermitian form on
$T_{\mathcal{A}_{W_2}/W_2}$.
Thus, since $\overline{W_2}$ is compact, there is a real number
$\lambda$ such that $0 < \lambda < 1$ and $\omega' - \lambda \omega$
is positive on $T_{\mathcal{A}_{W_2}/W_2}$.
To see our claim, clearly we may assume that
\[
W_2 = \DD^d = \{ (t_1, \ldots, t_d) \in \CC^d \mid |t_i| < 1 \}.
\] 
Let $\mu : \CC^g \times \DD^d \to \mathcal{A}_{\DD^d}$ 
be the universal covering
of $\mathcal{A}_{\DD^d}$ such that $\mu$ is a morphism over $\DD^d$, and
$\mu$ is a homomorphism on each fiber over $\DD^d$.
Let $(z_1, \cdots, z_g)$ be a coordinate of $\CC^g$.
Then, $\mu^*(\omega)$ and $\mu^*(\omega')$ can be written by the forms
\[
\sum_{ij} b_{ij}(z, t) d z_i \wedge d \bar{z}_j,
\]
where $b_{ij}(z, t)$'s are bounded $C^{\infty}$-functions.
Moreover, we set
\[
\begin{cases}
A_1 = \mu^*(c_1(\LL, \normabb_{cub})) - \mu^*(\omega) \\
A_2 = \mu^*(c_1(\LL, \normabb_0)) - \mu^*(\omega') \\
A_3 = \mu^*(c_1(\LL, \normabb_1)) -
(\pi \cdot \mu)^*(\phi) \mu^*(c_1(\LL, \normabb_{cub})) -
(1 - (\pi \cdot \mu)^*(\phi)) \mu^*(c_1(\LL, \normabb_0)).
\end{cases}
\]
Then, it is easy to see that each $A_i$'s can written by the form
\[
\sum_{ik} c_{ik}(z, t) d z_i \wedge d \bar{t}_k + 
\sum_{lj} c'_{lj}(z, t) d t_l \wedge d \bar{z}_j +
\sum_{lk} c''_{lk}(z, t) d t_l \wedge d \bar{t}_k,
\]
where $c_{ik}$'s, $c'_{lj}$'s and $c''_{lk}$'s are 
bounded $C^{\infty}$-functions.
Here let us see that $2^{-2n} [2^n]^*(A_i)$ ($i = 1, 2, 3$)
converges uniformly to $0$.
Indeed,
\begin{multline*}
2^{-2n} [2^n]^* \left(
\sum_{ik} c_{ik}(z, t) d z_i \wedge d \bar{t}_k + 
\sum_{lj} c'_{lj}(z, t) d t_l \wedge d \bar{z}_j +
\sum_{lk} c''_{lk}(z, t) d t_l \wedge d \bar{t}_k
\right)
= \\
\sum_{ik} 2^{-n}c_{ik}(2^n z, t) d z_i \wedge d \bar{t}_k + 
\sum_{lj} 2^{-n}c'_{lj}(2^n z, t) d t_l \wedge d \bar{z}_j +
\sum_{lk} 2^{-2n} c''_{lk}(2^n z, t) d t_l \wedge d \bar{t}_k
\end{multline*}
Thus, we have our assertion 
because $c_{ik}$'s, $c'_{lj}$'s and $c''_{lk}$'s are bounded.

Next, we try to see that $A_1 = 0$.
For, $2^{-2n}[2^n]^*(\mu^*(c_1(\LL, \normabb_{cub}))) = 
\mu^*(c_1(\LL, \normabb_{cub}))$ for all $n$, which shows us that
$2^{-2n}[2^n]^*(\omega) = \omega$ and
$2^{-2n}[2^n]^*(A_1) = A_1$. Hence, $A_1$ must be zero.

Thus, if we set $A = (1 - (\pi \cdot \mu)^*(\phi))A_2 + A_3$ and
\[
C =
(1-\lambda) (\pi \cdot \mu)^*(\phi)\mu^*(\omega)
+ (1 - (\pi \cdot \mu)^*(\phi))\mu^*(\omega' - \lambda \omega),
\]
then
\[
\mu^*(c_1(\LL, \normabb_1))  =
\lambda \mu^*(\omega) + A + C
\]
and $C$ is semipositive.
Hence,
\begin{multline*}
2^{-2n}[2^n]^*(\mu^*(c_1(\LL, \normabb_1))) + 
\epsilon (\pi \cdot \mu)^*(c_1(\overline{H}))
= \\
\lambda \mu^*(\omega) + \epsilon (\pi \cdot \mu)^*(c_1(\overline{H})) +
2^{-2n}[2^n]^*(A) + 2^{-2n}[2^n]^*(C).
\end{multline*}
On the other hand, $\lambda \mu^*(\omega) + 
\epsilon (\pi \cdot \mu)^*(c_1(\overline{H}))$
is positive and $2^{-2n}[2^n]^*(A)$ converges uniformly to $0$.
Thus,
\[
\lambda \mu^*(\omega) + \epsilon (\pi \cdot \mu)^*(c_1(\overline{H})) +
2^{-2n}[2^n]^*(A) 
\]
is positive if $n$ is sufficiently large.
Hence we get our claim.

\medskip
To get an invariant metric $\normabb$ under $F_{\infty}$,
over $F_{\infty}(W_1)$, we replace $\normabb_0$ by
$F_{\infty}^*(\normabb_1)$.
In this way, we have a hermitian line bundle 
$\overline{\LL} = (\LL, \normabb)$ with the following
properties:
\begin{enumerate}
\renewcommand{\labelenumi}{(\alph{enumi})}
\item
$c_1(\overline{\LL})$ is positive over 
$\mathcal{A}(\CC) \setminus \pi^{-1}(W_2 \cup F_{\infty}(W_2))$.

\item
For any $\epsilon > 0$, there is a positive number $n_0$ such that,
for all $n \geq n_0$,
\[
2^{-2n}([2]^n)^*(c_1(\overline{\LL})) + \epsilon \pi^* c_1(\overline{H})
\]
is positive on
$\pi^{-1}(W_2 \cup F_{\infty}(W_2))$.

\item
$2^{-2n} ([2]^n)^* c_1(\overline{\LL}) =  c_1(\overline{\LL})$ on $W$
for all $n$,
where $W = W_3$.
\end{enumerate}

\bigskip
Let $f_n : \mathcal{A}_n \to \mathcal{A}$ be the normalization of
\[
\mathcal{A}_U \overset{[2]^n}{\longrightarrow} \mathcal{A}_U
\hookrightarrow \mathcal{A}.
\]
Then, by projection formula,
\[
h^{\overline{B}}_{(\mathcal{A}_n, f_n^*(\overline{\LL}))}(x)
= h^{\overline{B}}_{(\mathcal{A}, \overline{\LL})}(2^n x)
\]
for all $x \in A(\overline{K})$.
Thus, if we set $\overline{\LL}'_{n} = 2^{-2n} f_n^*(\overline{\LL})$,
then,
\[
h^{\overline{B}}_{(\mathcal{A}_n, \overline{\LL}'_n)}(x)
= 2^{-2n} h^{\overline{B}}_{(\mathcal{A}, \overline{\LL})}(2^n x).
\]
Therefore,
\[
\lim_{n \to \infty} \sup_{x \in A(\overline{K})} 
\vert \hat{h}^{\overline{B}}_L(x) - 
h^{\overline{B}}_{(\mathcal{A}_n, \overline{\LL}'_n)}(x) \vert
= 0.
\]
Moreover,
\[
h^{\overline{B}}_{(\mathcal{A}_n, \overline{\LL}'_n + 
\epsilon \pi_n^*(\overline{H}))}
= h^{\overline{B}}_{(\mathcal{A}_n, \overline{\LL}'_n)}
\]
for any positive rational number $\epsilon$,
where $\pi_n : \mathcal{A}_n \to B$ is the canonical morphism.
Thus, if we set $\overline{\LL}_n =
\overline{\LL}'_n + \pi_n^*(\overline{H})$, then
a sequence of models $(\mathcal{A}_n, \LL_n)$ satisfies our
desired properties.
\QED

\section{Bogomolov's conjecture over finitely generated fields}
Let $K$ be a finitely generated field over $\QQ$
with $d = \trdeg_{\QQ}(K)$, 
and $\overline{B} = (B; \overline{H}_1, \ldots, \overline{H}_d)$
a polarization of $K$.
Let $A$ be an abelian variety over $K$.
In this section, we would like to prove the following theorem,
which is a generalization of
results due to Ullmo \cite{UlPos} and Zhang \cite{ZhEqui}.

\begin{Theorem}
\label{thm:bogomolov:conj:fun}
Let $X$ be a subvariety of $A_{\overline{K}}$, and
$L$ a symmetric ample line bundle on $A$.
We assume that $\overline{B}$ is big, i.e.,
$\overline{H}_i$'s are nef and big.
If the set 
\[
 \{ x \in X(\overline{K}) \mid \hat{h}_L^{\overline{B}}(x)
   \leq \epsilon \}
\]
is Zariski dense in $X$ for any $\epsilon > 0$,
then $X$ is a translation of an abelian subvariety of $A_{\overline{K}}$
by a torsion point.
\end{Theorem}

\Proof
First of all, note that in order to prove our theorem,
we can replace the field $K$ by a finite extension of $K$
if it is necessary.

We set
\[
G(X) = \{ a \in A(\overline{K}) \mid a + X = X \}.
\]
First, let us consider the case where $G(X)$ is trivial.
In this case, we need to show that
there is a torsion point $x$ of $A$ with $X = \{ x \}$.
For this purpose, it is sufficient to show that $\dim X = 0$.
For, if we set $X = \{ x \}$, then
$\hat{h}_L^{\overline{B}}(x) = 0$.
Thus, $x$ is a torsion point by
Proposition~\ref{prop:positivity:canonical:height}. 

From now on, we assume that $\dim X \geq 1$.
Changing $K$ by a finite extension of $K$, if necessarily,
by Proposition~\ref{prop:exist:line:pencil:type},
we may assume that there is a $C^{\infty}$-hermitian line bundle $\overline{H}_0$
with $\adeg(\acherncl_1(\overline{H}_0)^{d+1}) = 0$ and
$\deg((H_0)_{\QQ}^d) > 0$.
Let $\overline{B}_0 = (B; \overline{H}_0, \ldots, \overline{H}_0)$ be
a polarization of $K$ given by $\overline{H}_0$.
Then, by virtue of (3) of Proposition~\ref{prop:comp:canonical:height},
there is a positive constant $a$ with 
$\hat{h}^{\overline{B}_0}_L \leq a \hat{h}^{\overline{B}}_L$.
Thus, the set
\[
 \{ x \in X(\overline{K}) \mid \hat{h}^{\overline{B}_0}_L(x) \leq
\epsilon \}
\]
is Zariski dense for any $\epsilon > 0$.
Therefore, we will try to find a contradiction using hypotheses:
\begin{enumerate}
\renewcommand{\labelenumi}{(\alph{enumi})}
\item
$G(X) = \{ 0 \}$.

\item
$\dim X \geq 1$.

\item
The set 
$ \{ x \in X(\overline{K}) \mid \hat{h}^{\overline{B}_0}_L(x) \leq
\epsilon \}$
is Zariski dense for any $\epsilon > 0$.
\end{enumerate}

Here we consider a morphism
\[
\phi_m : A_{\overline{K}}^m \to A_{\overline{K}}^{m-1}
\]
given by 
$\phi_m(x_1, \ldots, x_m) = (x_1 - x_2, x_2 - x_3, \ldots, x_{m-1} - x_m)$.
Then, since $G(X) = \{ 0 \}$, in the same way
as the proof of \cite[Lemma~3.1]{ZhEqui}, we can see that
if $m$ is sufficiently large,
then $\phi_m$ induces
a birational morphism $X^m \to \phi_m(X^m)$.
Considering a finite extension of $K$, we may assume that
$X$ is defined over $K$, and that
$\phi_m$ induces a birational morphism $X^m \to \phi_m(X^m)$ over $K$.
Here, if it is necessary, we change
$\overline{B}_0$ by the polarization induced by
the extension of $K$ accordingly.

We note that the above hypothesis (c) does not depend on the choice of 
the ample and symmetric line bundle $L$ 
by virtue of (2) of Proposition~\ref{prop:comp:canonical:height}.
Hence, by Proposition~\ref{prop:construction:good:seq:model},
there is a sequence of $C^{\infty}$-models $(\mathcal{A}_n, \overline{\LL}_n)$
of $(A, L)$ with the following properties.
\begin{enumerate}
\renewcommand{\labelenumi}{(\arabic{enumi})}
\item
There is a Zariski open set $U$ of $B$ such that
$(\mathcal{A}_n)_U = (\mathcal{A}_1)_U$ for all $n$, and that
$(\mathcal{A}_1)_U \to U$ is an abelian scheme over $U$.

\item
If $n$ is sufficiently large, then
$\LL_n$ is ample and
$\overline{\LL}_n$ is vertically nef.

\item
${\displaystyle \lim_{n \to \infty} 
\sup_{x \in A(\overline{K})} 
\vert \hat{h}^{\overline{B}_0}_L(x) - 
h^{\overline{B}_0}_{(\mathcal{A}_n, \overline{\LL}_n)}(x) \vert
= 0}$.

\item
There are a connected open set $W$ of $U(\CC)$,
and a positive $C^{\infty}$-form $\omega$ on $(\mathcal{A}_1)_W$ such that
$W$ is non-singular, $c_1(\overline{H}_0)$ is positive on $W$,
and that $c_1(\overline{\LL}_n) = \omega$ on 
$(\mathcal{A}_1)_W$ for all $n \gg 0$.
\end{enumerate}

For simplicity, we denote $\mathcal{A}_1$ by $\mathcal{A}$.
Let $\mathcal{A}^m$ (resp. $\mathcal{A}^m_n$) be the main component of
\[
\underbrace{\mathcal{A} \times_B \cdots \times_B 
\mathcal{A}}_{\text{$m$-times}}
\qquad
\left( \text{resp.}\quad
\underbrace{\mathcal{A}_n \times_B \cdots \times_B 
\mathcal{A}_n}_{\text{$m$-times}}
\right).
\]
Let $\pi^m : \mathcal{A}^m \to B$ and
$\pi^m_n : \mathcal{A}^m_n \to B$ be the canonical projections.
Let $\XX^m$ (resp. $\XX_n^m$) be the closure
of $X^m$ in $\mathcal{A}^m$ (resp. $\mathcal{A}^m_n$).
Further, let $\mathcal{Y}^m$ (resp. $\mathcal{Y}^m_n$) 
be the closure of $\phi_m(X^m)$
in $\mathcal{A}^{m-1}$ (resp $\mathcal{A}^{m-1}_n$).
We set
\[
\overline{\mathcal{P}}^m = p_1^*(\overline{\LL}) \otimes \cdots 
\otimes p_m^*(\overline{\LL})
\quad\text{and}\quad
\overline{\mathcal{P}}^m_n = p_1^*(\overline{\LL}_n) \otimes \cdots 
\otimes p_m^*(\overline{\LL}_n),
\]
where
$p_i$ is the projection to the $i$-th factor.
Note that $\phi_m$ extends to $\mathcal{A}_U^m \to \mathcal{A}_U^{m-1}$.
By abuse of notation, we denote this extension by $\phi_m$.
Then, $\phi_m$ induces a birational morphism
$\XX_U^m \to \mathcal{Y}^m_U$.
Let $V$ and $V'$ be Zariski open sets of $\XX^m_U$ and $\mathcal{Y}^m_U$
respectively such that $\phi_m$ gives rise to an isomorphism
$V \overset{\sim}{\longrightarrow} V'$.

Since $X$ has only countably many subvarieties over $\overline{K}$,
let $\{ Y_t \}_{t=1}^{\infty}$ be the set of all
proper subvarieties of $X$. By the hypothesis (c),
we can find $x_t \in X(\overline{K})$ such that
$x_t \not\in \bigcup_{i=1}^t Y_i$ and 
$\hat{h}^{\overline{B}_0}_L(x_t) \leq 1/t$.
Then, we have a generic sequence
$\{ x_t \}$ of $X(\overline{K})$
with ${\displaystyle \lim_{t \to \infty} \hat{h}^{\overline{B}_0}_L(x_t) = 0}$.
Let us fix a bijection $\tau : \NN \to \NN^m$.
We denote $(x_{\tau_1(t)}, \ldots, x_{\tau_m(t)}) \in X^m$
by $x_{\tau(t)}$, where $\tau(t) = (\tau_1(t), \ldots, \tau_m(t))$.
Since $\{ x_{\tau(t)} \}$ is Zariski dense in $X^m$,
in the same way as before, we can find a generic subsequence
of $\{ x_{\tau(t)} \}$. Thus, we may assume
that $\{ x_{\tau(t)} \}$ is a generic sequence.
Moreover, considering a subsequence of $\{ x_t \}$,
we may further assume that $x_{\tau(t)} \in V_{K}$.
Further, we can see that $\lim_{t \to \infty} 
\hat{h}^{\overline{B}_0}_{P^m}(x_{\tau(t)}) = 0$ and
$\hat{h}^{\overline{B}_0}_{P^{m-1}}(\phi_m(x_{\tau(t)})) = 0$,
where $P^m = \mathcal{P}^m_K$ and $P^{m-1} = \mathcal{P}^{m-1}_K$.
Thus, using the equidistribution theorem
(cf. Theorem~\ref{thm:equi:dist}), 
over $(\pi^m)^{-1}(W) \cap \XX^m$,
\[
\lim_{t\to\infty}
\frac{\delta_{\Delta_{x_{\tau(t)}}} \wedge (\pi^m)^*(c_1(\overline{H}_0)^{\wedge d})}
{\deg(\Delta_{x_{\tau(t)}} \to W)} 
=
\left[ \frac{(p_1^*(\omega) + \cdots + p_m^*(\omega))^{\wedge em} \wedge 
(\pi^m)^*(c_1(\overline{H}_0)^{\wedge d})}
{\deg\left( \left(\rest{P^m_K}{X^m}\right)^{em} \right)} \right],
\]
where $e = \dim X$.
Moreover, if we denote $\alpha$ by the restriction of
$p_1^*(\omega) + \cdots + p_{m-1}^*(\omega)$ to 
$(\pi_{m-1})^{-1}(W) \cap \mathcal{Y}^m$,
then
\[
\lim_{t\to\infty}
\frac{\delta_{\Delta_{\phi_m(x_{\tau(t)})}} \wedge 
(\pi^{m-1})^*(c_1(\overline{H}_0))^{\wedge d}}
{\deg(\Delta_{\phi_m(x_{\tau(t)})} \to W)} 
=
\left[
\frac{\alpha^{\wedge em} \wedge (\pi^{m-1})^*(c_1(\overline{H}_0)^{\wedge d})}
{\deg \left( \left( \rest{P^{m-1}_K}{Y^m}\right)^{em} \right)}
\right],
\]
where $Y^m = \mathcal{Y}^m_K$. Note that
on $V \cap (\pi^m)^{-1}(W)$,
\[
\frac{\delta_{\Delta_{x_{\tau(t)}}} \wedge (\pi^m)^*(c_1(\overline{H}_0)^{\wedge d})}
{\deg(\Delta_{x_{\tau(t)}} \to W)}
\quad\text{and}\quad
\frac{\delta_{\Delta_{\phi_m(x_{\tau(t)})}} \wedge 
(\pi^{m-1})^*(c_1(\overline{H}_0))^{\wedge d}}
{\deg(\Delta_{\phi_m(x_{\tau(t)})} \to W)}
\]
give rise to the same current via
the isomorphism $\phi_m : V \overset{\sim}{\longrightarrow} V'$.
Thus two limit $C^{\infty}$-forms 
\[
\frac{(p_1^*(\omega) + \cdots + p_m^*(\omega))^{\wedge em} \wedge 
(\pi^m)^*(c_1(\overline{H}_0)^{\wedge d})}
{\deg\left( \left(\rest{P^m_K}{X^m}\right)^{em} \right)}
\quad\text{and}\quad
\phi_m^* \left(
\frac{\alpha^{\wedge em} \wedge (\pi^{m-1})^*(c_1(\overline{H}_0)^{\wedge d})}
{\deg \left( \left( \rest{P^{m-1}_K}{Y^m}\right)^{em} \right)}
\right)
\]
are same forms on $V \cap (\pi^m)^{-1}(W)$.
Hence we have
\[
(p_1^*(\omega) + \cdots + p_m^*(\omega))^{\wedge em} \wedge 
(\pi^m)^*(c_1(\overline{H}_0)^{\wedge d}) =
c (\phi_m^*(\alpha^{\wedge em})) \wedge (\pi^m)^*(c_1(\overline{H}_0)^{\wedge d})
\]
over $(\pi^{m})^{-1}(W)\cap \XX^m$ for some positive constant $c$.
Let $\mathcal{X}$ be the closure of $X$ in $\mathcal{A}$.
We choose a general point $w$ of $W$.
Then, we can easily see that
\[
\rest{(p_1^*(\omega) + \cdots + p_m^*(\omega))^{\wedge em}}{(\XX_w)^m} =
\rest{c (\phi_m^*(\alpha^{\wedge em}))}{(\XX_w)^m}
\]
because $c_1(\overline{H}_0)$ is positive on $W$,
where $\mathcal{X}_w$ is the fiber of
$\mathcal{X} \to B$ over $w$.
This is a contradiction because the left hand side is positive, but
the right hand side is not positive along the diagonal of
$(\XX_w)^m$.

\bigskip
Next let us consider a general case.
Let $\mu : A_{\overline{K}} \to A' = A_{\overline{K}}/G(X)$
be the quotient of $A_{\overline{K}}$ by $G(X)$.
We set $X' = \mu(X)$. Then, it is easy to see that
$\mu^{-1}(X') = X$ and $G(X') = \{ 0 \}$.
Let $L'$ be a symmetric and ample line bundle on $A'$.
Let $K'$ be a finite extension field of $K$
such that $\mu$, $A'$, and $L'$ are defined over $K'$.
Let $\overline{B}'$ be the polarization of $K'$ induced by
$\overline{B}$.
Since $L$ is ample, 
by Proposition~\ref{prop:base:change:field} and 
(2) of Proposition~\ref{prop:comp:canonical:height},
there is a positive number $a$ such that
$\hat{h}^{\overline{B}'}_{\mu^*(L')} \leq a \hat{h}^{\overline{B}}_L$.
Thus,
\begin{align*}
\{ x \in X(\overline{K}) \mid 
\hat{h}^{\overline{B}}_{L}(x) \leq \epsilon/a \}
& \subseteq
\{ x \in X(\overline{K}) \mid 
\hat{h}^{\overline{B}'}_{\mu^{*}(L')}(x) \leq \epsilon \} \\
& =
\mu^{-1}\left(
\{ x' \in X'(\overline{K}) \mid \hat{h}_{L'}^{\overline{B}'}(x')
\leq \epsilon \}
\right).
\end{align*}
Therefore, the set
$\{ x' \in X'(\overline{K}) \mid \hat{h}_{L'}^{\overline{B}'}(x')
\leq \epsilon \}$
is Zariski dense in $X'$ for any $\epsilon > 0$.
Thus, by the previous observation,
$X' = \{ x' \}$ for some torsion point $x'$ of $A'$.
Hence, $X$ is a coset of $G(X)$ because $\mu^{-1}(X') = X$.
In particular, $G(X)$ is an abelian subvariety.
Thus, it is sufficient to show that
there is a torsion point $x$ of $A$ with $\mu(x) = x'$.
First, pick up $x_1$ of $A$ with $\mu(x_1) = x'$.
Since $x'$ is a torsion point, there is a positive number $n$
with $n x_1 \in G(X)$.
Here $G(X)$ is a divisible group. Thus, we can find $x_2 \in G(X)$
with $n x_1 = n x_2$. Hence, if we set $x = x_1 - x_2$,
then we have a desired torsion point.
\QED

As corollary, we can recover the following Raynaud's result
(\cite{Ray1} and \cite{Ray2}).

%%
%\begin{Corollary}
%\label{cor:raynaud:thm}
%Let $F$ be an algebraically closed field of characteristic zero.
%Let $A$ be an abelian variety over $F$, and
%$X$ a subvariety of $A$. If $X(F) \cap A(F)_{tor}$ is Zariski dense
%in $X(F)$, then $X$ is a translation of an abelian subvariety of $A$
%by a torsion point.
%\end{Corollary}
%%
\begin{Corollary}
\label{cor:raynaud:thm}
Let $A$ be an abelian variety over 
an algebraically closed field $F$ of characteristic zero, and
$Z$ a reduced subscheme of $A$. 
Then, every irreducible component of
the Zariski closure of $Z(F) \cap A(F)_{tor}$ in A is
a translation of an abelian subvariety of $A$
by a torsion point.
Consequently, there are finitely many abelian subvarieties
$B_1, \ldots, B_n$ of $A$ and torsion points $b_1, \ldots, b_n$ of $A(F)$
such that
\[
\overline{Z(F) \cap A(F)_{tor}} = \bigcup_{i=1}^n (B_i(F) + b_i)
\quad\text{and}\quad
Z(F) \cap A(F)_{tor} = \bigcup_{i=1}^n (B_i(F)_{tor} + b_i).
\]
\end{Corollary}

%%
%\Proof
%Let $K$ be a subfield of $F$ such that
%$K$ is a finitely generated field over $\QQ$, and that
%$A$ and $X$ are defined over $K$.
%Then, $A(F)_{tor} = A(\overline{K})_{tor}$.
%Thus, we can see that $X(\overline{K}) \cap A(\overline{K})_{tor}$
%is Zariski dense in $X(\overline{K})$.
%Therefore, by Theorem~\ref{thm:bogomolov:conj:fun},
%$X$ is a translation of an abelian subvariety of $A$
%by a torsion point. 
%\QED
%%

\Proof
Let $X$ be an irreducible component of
the Zariski closure of $Z(F) \cap A(F)_{tor}$ in A.
Then, it is easy to see that $X(F) \cap A(F)_{tor}$
is Zariski dense in $X$.
Let $K$ be a subfield of $F$ such that
$K$ is a finitely generated field over $\QQ$, and that
$A$ and $X$ are defined over $K$.
Then, $A(F)_{tor} = A(\overline{K})_{tor}$.
Thus, we can see that $X(\overline{K}) \cap A(\overline{K})_{tor}$
is Zariski dense in $X(\overline{K})$.
Therefore, by Theorem~\ref{thm:bogomolov:conj:fun},
$X$ is a translation of an abelian subvariety of $A$
by a torsion point. 
\QED

\begin{Remark}
We assume that $\trdeg_{\QQ}(K) = 1$.
Then, as in the introduction, we have two types of heights
$\hat{h}^{geom}_L$ and $\hat{h}^{arith}_L$.
Then, by virtue of (3) of Proposition~\ref{prop:comp:canonical:height},
there is a positive constant $a$
with $\hat{h}^{geom}_L \leq a \hat{h}^{arith}_L$.
This means that Bogomolov's conjecture for the geometric height
$\hat{h}^{geom}_L$ is a subtle problem. Actually, if
the trace of $A$ is not zero, the conjecture does not holds in general.
(For example, take $X$ as a non-torsion point $P$ with $\hat{h}^{geom}_L(P) = 0$.)
However, we can expect the conjecture for $\hat{h}^{geom}_L$ if
either the trace of $A$ is zero, or
$X$ is a non-isotrivial curve and $A$ is its jacobian.
For example, in \cite{MoRB}, the author gives an affirmative answer under
the assumption of singular fibers on the stable model of $X$.
\end{Remark}


\bigskip
\begin{thebibliography}{99}

\bibitem{Alt}
A. Altman, The size function on abelian varieties,
Trans. of A.M.S., 164(1972), 153--161.

\bibitem{GSArInt}
H. Gillet and C. Soul\'{e},
Arithmetic Intersection Theory,
Publ. Math. (IHES), 72 (1990), 93--174.

\bibitem{Hiro}
H. Hironaka,
Resolution of singularities of an algebraic variety over
a field of characteristic zero,
Ann. of Math. 79 (1964), 109--326.

\bibitem{KMSemi}
S. Kawaguchi and A. Moriwaki,
Inequalities for semistable families for arithmetic varieties,
(alg-geom/9710007).

\bibitem{LaFund}
S. Lang,
Fundamentals of diophantine geometry,
(1983), Springer.

\bibitem{MoBMet}
L. Moret-Bailly,
M\'{e}triques permises,
S\'{e}minaire sur les pinceaux arithm\'{e}tiques:
La Conjecture de Mordell,
Ast\'{e}risque 127 (1985), 29-87.

\bibitem{MoRB}
A. Moriwaki,
Relative Bogomolov's inequality and the cone of positive divisors
on the moduli space of stable curves,
J. of AMS, 11 (1998), 569--600.

\bibitem{Ray1}
M. Raynaud,
Courbes sur une vari\'{e}t\'{e} ab\'{e}lienne et points de torsion,
Invent. math., 71 (1983), 207--233.

\bibitem{Ray2}
M. Raynaud,
Sous-vari\'{e}t\'{e} d'une vari\'{e}t\'{e} ab\'{e}lienne et
points de torsion,
in Arithmetic and Geometry vol.1, (1983).

\bibitem{SerPG}
J. -P. Serre,
Propri\'{e}t\'{e}s galoisiennes des points
d'ordre fini des courbes elliptiques,
Invent. math., 15 (1972), 259--331.

\bibitem{SoAr}
C. Soul\'{e} et al,
Lectures on Arakelov Geometry,
Cambridge studies in advanced mathematics, 33,
Cambridge University Press.

\bibitem{SUZEqui}
L. Szpiro, E. Ullmo, and S. Zhang,
Equir\'{e}partition des petits points,
Invent. math., 127 (1997), 337-347.

\bibitem{UlPos}
E. Ullmo,
Positivit\'{e} et disc\'{e}tion des points alg\'{e}briques
des courbes,
Ann. Math., 147 (1998), 167-179.

\bibitem{ZhPos}
S. Zhang,
Positive line bundles on arithmetic varieties,
J. of AMS., 8 (1995), 187--221.

\bibitem{ZhEqui}
S. Zhang,
Equidistribution of small points on abelian varieties,
Ann. Math., 147 (1998), 159-165.
\end{thebibliography}
\end{document}